\documentclass[a4paper,12pt]{article}
\usepackage[T1]{fontenc}
\usepackage[utf8]{inputenc}
\usepackage[english]{babel}
\usepackage{amsmath,amssymb}
\usepackage{amsthm,amsrefs}
\usepackage[autolanguage,np]{numprint}
\usepackage{graphicx}
\usepackage{subfig,float}

\newtheorem{lemma}{Lemma}
\newtheorem{proposition}{Proposition}
\newtheorem{theorem}{Theorem}
\newtheorem{corollary}{Corollary}
\theoremstyle{remark}
\newtheorem{remark}{Remark}
\newtheorem{definition}{Definition}
\DeclareMathOperator{\Dom}{Dom}
\DeclareMathOperator{\sgn}{sgn}
\DeclareMathOperator{\grandO}{O}
\DeclareMathOperator{\argmax}{argmax}
\DeclareMathOperator{\Var}{Var}
\newcommand{\RR}{\mathbb{R}}
\newcommand{\PP}{\mathbb{P}}
\newcommand{\EE}{\mathbb{E}}
\newcommand{\NN}{\mathbb{N}}
\newcommand{\vd}{\,\mathrm{d}}
\newcommand{\dd}{\mathrm{d}}
\newcommand{\cN}{\mathcal{N}}
\newcommand{\cC}{\mathcal{C}}
\newcommand{\indicator}{\mathbf{1}}

\title{Is a Brownian motion skew?}

\author{
Anotine Lejay\thanks{Projet TOSCA (INRIA and Institut \'Elie Cartan, UMR CNRS,
Nancy-Universit\'e, CNRS, INRIA), IECN, BP 239, 54506 Vand\oe
uvre-l\`es-Nancy \textsc{cedex}, France;
email: $\langle$\texttt{Antoine.Lejay@iecn.u-nancy.fr}$\rangle$.
This author  has been supported by the MathAmSud program.}
\and 
Ernesto Mordecki\thanks{Centro de Matem\'atica, Facultad de Ciencias,
Universidad de la Rep\'ublica. Igu\'a 4225, 11400, Montevideo,
Uruguay; email: $\langle$\texttt{mordecki@cmat.edu.uy}$\rangle$.
The author's research is supported by CMAT-UDELAR, and by project \emph{SAMP Network}, of the MathAmSud program.} 
\and 
Soledad Torres\thanks{CIMFAV-DEUV; 
Universidad de Valpara\'{\i}so; Casilla 123-V; Chile; 
email: $\langle$\texttt{soledad.torres@uv.cl}$\rangle$. The author's research is supported by 
PBCT-ACT 13 Stochastic Analysis Laboratory, Chile.}
}
\date{\today}

\begin{document}

\maketitle

\begin{abstract}
We study the asymptotic behavior of the maximum likelihood estimator corresponding to the observation of a trajectory of a Skew Brownian motion, through a uniform time discretization.  We characterize the speed of convergence and the limiting distribution when the step size goes to zero, which in this case are non-classical, under the null hypothesis of the Skew Brownian motion being an usual Brownian motion.  This allows to design a test on the skewness parameter.  We show that numerical simulations that can be easily performed to estimate the skewness parameter, and provide an application in Biology.
\end{abstract}

\noindent\textbf{Keywords: } Skew Brownian motion, statistical estimation, 
maximum likelihood.

\section{Introduction}


The Skew Brownian Motion (SBm) has attracted interest within other facts,
due to its relations with diffusions with discontinuous coefficients or to media with permeable barriers,
being the first example of the solution of a stochastic differential equation with the local time of the solution as drift \cite{harrison}: 
the SBm 
$X=\{X_t\colon 0\leq t\leq T\}$ 
can be defined as the strong solution of the stochastic differential equation
\begin{equation}
\label{eq:skb}
X_t=x+B_t+\theta \ell^x_t,
\end{equation}
where $B=\{B_t\colon 0\leq t\leq T\}$ is a standard Brownian motion defined on a probability space 
$(\Omega,\mathcal{F},\PP)$,
the initial condition is $x\geq0 $ (the case $x<0$ is symmetrical), 
$\theta\in[-1,1]$ is the  \emph{skewness parameter},
and $\ell^x=\{\ell^x_t\colon 0\leq t\leq T\}$ is the
local time at level zero of the (unknown) 
solution $X$ of the equation departing from $x$, defined by
\begin{equation}
\label{eq:local}
\ell^x_t = \lim_{\epsilon \to 0}{1 \over 2 \epsilon} \int_0^t 1_{(-\epsilon,
\epsilon)}(X_s)ds.
\end{equation}
In case $x=0$ we denote $\ell^0_t=\ell_t$, being this case particularly interesting due to some explicit calculations 
that can be carried out (see Section \ref{section:ld}).


In the literature the skewness parameter is sometimes defined as $p=(\theta+1)/2$;
this second parametrization being more convenient for an alternative construction of the SBm:
depart from the reflected Brownian motion and choose, 
independently with probability $p\in[0,1]$, 
whether each particular excursion of the reflected Brownian motion remains positive.


In the special case $\theta = 1\  (p=1)$, 
the solution to \eqref{eq:skb} is the reflected Brownian motion. 
The case $\theta=0\ (p=1/2)$ corresponds to the the standard Brownian motion.

Recently, several papers have considered the SBm in modelling or simulation issues, 
as well as some optimization problems.
See the review by A.~Lejay \cite{lejay} for references on the subject, 
as well as a survey of the various possible constructions and applications of the SBm.


In this paper we are interested in the statistical estimation of $\theta$, the skewness parameter,
when we observe a trajectory of the process through an equally spaced time grid. 
From the statistical point of view we find this problem interesting because it is
in certain sense intermediate between the classical problem of drift estimation in a diffusion, 
where the measures generated by the trajectories of the process for different values of the parameter are equivalent
\cites{K,lipster}, 
and the estimation of the variance (the volatility in financial terms) of a diffusion 
(see for instance \cite{florens1}, or \cite{jacod2} and the
references therein), 
where the probability measures generated by the trajectories are singular for different values of the parameter.
At the best of our knowledge, the only estimator
of $\theta$ is the one constructed by O.~Bardou
and M.~Martinez~\cite{bardou}, where they assume
that the SBm is reflected at levels $1$ and $-1$ to ensure ergodicity, considering 
a different scheme of observation of the trajectory.


Our main result states that the maximum likelihood estimator
(MLE) corresponding to the observation of a discretization
of one trajectory of the process, with the corresponding
normalization, satisfies
the so called Local Asymptotically Mixed Normality (LAMN) property at the point $\theta=0$. With this result and the identification of the limiting distribution 
of the scaled MLE estimator, one may construct some hypothesis test to determine
whether or not the Brownian is skew.
This fact suggests certain asymptotic properties
of the MLE, as exposed for instance in the classical book of Ibragimov and Has'minskii \cite{ih}.
Nevertheless, 
as our results in terms of convergence of statistical experiments are not exactly the ones
needed in the hypothesis of general LAMN theorems,  
we follow a direct approach to construct the estimator and to study its asymptotic properties. 
This approach, that can be followed in rare occasions, has the
advantage of clarifying the proof of the asymptotic properties 
and providing insight in the corresponding numerical computations.


The rest of the paper is organized as follows. Section 2 describes
the maximum likelihood methodology and the convergence results. 
In Section 3 we describe the limit distribution. Sections 4 presents 
the statistical Test and some numerical simulations on the likelihood function.
Section 5 presents an application to diffusion of species in two different habitats, 
and Section 6 our conclusions. 
Finally, in the Appendix we provide the 
theorems taken from \cite{jacod} used in the proof of our main results in Section 2.

\section{The maximum likelihood estimator}\label{section:mle}

Consider the SBm $X$ with parameter $\theta\in(-1,1)$ defined in \eqref{eq:skb} 
and the sampling scheme 
denoted by $X_i:=X_{iT/n}\ (i=0,\dots,n)$, and $\Delta=T/n$.
In this section we derive the asymptotic behaviour of maximum
likelihood estimator $\theta_n$ of the 
parameter $\theta$ when we observe the sample $X_1, \ldots , X_n$. 
The transition density of the SBm of parameter $\theta\in[-1,1]$
is given by: 
\begin{equation*}
q_{\theta}(t,x,y)=p(t,y-x)+\sgn(y)\theta p(t,|x|+|y|),
\end{equation*}
where
\begin{equation*}
p(t,x)=\frac{1}{\sqrt{2\pi t}}\exp\Big(-\frac {x^2}{2t}\Big)
\end{equation*}
is the density of a Gaussian random variable with variance $t$ and mean $0$.
The likelihood of the sample is given by
\begin{align*}
Z_n(\theta)&=\prod_{i=0}^{n-1}{q_\theta(\Delta,X_i,X_{i+1})\over q_0(\Delta,X_i,X_{i+1})}.
\end{align*}
Observing that for any $x,y\in\RR$ we have 
\begin{align}
\frac{p(\Delta,|x|+|y|)}{q_0(\Delta,x,y)}
&=\exp\left(-\frac{|xy|+xy}{\Delta}\right)=\exp\left(\frac{-2(xy)^+}{\Delta}\right)\leq 1,
\label{eq:bound}
\end{align}
(where $z^+=(|z|+z)/2=\max(z,0)$), we can write
\begin{align*}
Z_n(\theta)&=
\prod_{\substack{X_i>0\\X_{i+1}<0}}(1-\theta)
\prod_{\substack{X_i<0\\X_{i+1}>0}}(1+\theta)
\prod_{\substack{X_i<0\\X_{i+1}<0}}\left(1-\theta e^{-2X_iX_{i+1}\over\Delta}\right)
\\ &\quad\times
\prod_{\substack{X_i>0\\X_{i+1}>0}}\left(1+\theta e^{-2X_iX_{i+1}\over\Delta}\right)
=
\prod_{i=0}^{n-1}\left(1+h(\sqrt{n}X_i,\sqrt{n}X_{i+1})\right),
\end{align*}
where
\begin{equation*}
h(x,y)=\sgn(x+y)\exp\left(-(2/T)(x(x+y))^+\right),
\end{equation*}
to see that $Z_n(\theta)$ is a polynomial of degree $n$, with $n$ real roots. 
Remember that we assume $X_0=x\geq 0$. In case the trajectory we observe
does not hit the zero level, we obtain
\begin{align*}
Z_n(\theta)&=
\prod_{i=0}^{n-1}\left(1+\theta e^{-2X_iX_{i+1}\over\Delta}\right)
\end{align*}
and $Z_n(\theta)$ is increasing in $\theta$. In this case our maximum likelihood 
estimator is $\theta_n=1$.  
In the case that the trajectory crosses the zero level, we see that
the polynomial has roots $\theta=\pm 1$ (for large enough $n$), and no roots
inside this interval. As $Z_n(0)=1$, this gives a unique maximum at the point $\theta_n$ 
in the interval $(-1,1)$. 

Our main result is the weak convergence of the MLE to a distribution that
we characterize. Three main differences can be noted in respect to the classical
statistical situation: (i) the convergence of the estimator is more slowly ($n^{1/4}$)
than in the classical case; (ii) the limit is not Gaussian, but a mixture of Gaussian random
variables; (iii) the convergence is stable, stronger than the usual convergence in distribution,
but natural in this context, known as local asymptotic mixed normality (LAMN) in the literature 
(see for instance \cite{lecam-yang}).
We also have to take into account, in accordance to our
previous discussions on the existence and the value of the MLE, that the event that the
trajectory hits the level zero is crucial in the results we obtain (in fact, if the 
trajectory does not hit this level, the MLE remains constant for all $n$). 
Consider then the events
\begin{equation}\label{eq:events}
A_n=\{\omega\colon\inf_{1\leq i\leq n}X_i<0\},
\qquad
A=\{\omega\colon\inf_{0\leq t\leq T}X_t(\omega)<0\}.
\end{equation}
As $X$ is continuous, $\indicator_{A_n}\to\indicator_A$ a.s. 
($\indicator_B$ stands for the indicator of the set $B$).
We now review the stable convergence (see \cite{jacod-shiryaev:1987}), 
and introduce the conditional stable convergence that will take place in our case.
\begin{definition} Consider a sequence of random variables 
$Y,Y_1,Y_2,\dots$ defined on a the probability space 
$(\Omega,\mathcal{F},\PP)$, and a $\sigma$-algebra 
$\mathcal{G}\subset\mathcal{F}$.

We say that the sequence of random variables 
$Y_1,Y_2\dots$ converge $\mathcal{G}$-stably in distribution to $Y$, and denote
\begin{equation*}
Y_n\xrightarrow[n\to\infty]{\text{$\mathcal{G}$-stably}}Y 
\end{equation*}
when
\begin{equation*}
\EE\left(Z f(Y_n)\right)\xrightarrow[n\to\infty]{}\EE\left(Z f(Y)\right)
\end{equation*}
for any bounded $\mathcal{G}$ measurable random variable $Z$, 
and any bounded and continuous function $f$.

Furthermore, consider a sequence of sets $A,A_1,A_2,\dots$.
We say that the sequence of random variables 
$Y_1,Y_2\dots$ conditional on $A_n$ converge 
$\mathcal{G}$-stably
in distribution to $Y$ conditional on $A$, and denote
\begin{equation*}
Y_n\mid A_n\xrightarrow[n\to\infty]{\text{$\mathcal{G}$-stably}}Y\mid A, 
\end{equation*}
when
\begin{equation*}
\EE\left(Z f(Y_n)\mid A_n\right)\xrightarrow[n\to\infty]{}\EE\left(Z f(Y)\mid A\right)
\end{equation*}
for any bounded $\mathcal{G}$ measurable random variable $Z$, 
and any bounded and continuous function $f$.
\end{definition}

We are now in position of presenting our main result. Indeed, 
this theorem will be an immediate sequel of Theorem~\ref{thm:2} below.

\begin{theorem}
\label{thm:1}
Consider a Skew Brownian motion defined in \eqref{eq:skb}
with the sampling scheme described in the beginning of Section \ref{section:mle} 
and the events $A_n$ and $A$ defined in \eqref{eq:events}. Then for 
the maximum likelihood estimator $\theta_n$ we have the convergence
\begin{equation*}
n^{1/4}\theta_n\mid A_n\xrightarrow[n\to\infty]{\text{$\mathcal{F}$-stably}}
{W(\ell^x_T)\over\ell^x_T}\mid A, 
\end{equation*}
under the Brownian motion distribution (that is when $\theta=0$), 
where $W=\{W_t\colon t\geq 0\}$ is a standard Brownian motion independent of $B$.
In particular, when $x=0$, we have
\begin{equation}
\label{eq:theta}
n^{1/4}\theta_n\xrightarrow[n\to\infty]{\text{$\mathcal{F}$-stably}}
{W(\ell_T^x)\over\ell_T^x} 
\end{equation}
\end{theorem}

\subsection{Some results on derivatives of the log-likelihood}

In order to study the asymptotic behaviour of  $\theta_n$, the MLE, 
we consider the \emph{log-likelihood}, defined by 
\begin{equation}
\label{eq:logLike}
L_n(\theta)=\log \prod_{i=0}^{n-1} q_\theta(\Delta,X_i,X_{i+1})
\end{equation}
and introduce its scaled (for notational convenience) $k$-th derivatives, for $k\geq 1$, by
\begin{equation*}
L_n^{(k)}(\theta)={1\over (k-1)!}{\partial^k\over\partial\theta^{k}}L_n(\theta),
\end{equation*}
that are computed as
\begin{equation}\label{derk}
L_n^{(k)}(\theta)=(-1)^{k-1}\sum_{i=0}^{n-1} 
\frac{\sgn(X_{i+1})^k p(\Delta,|X_i|+|X_{i+1}|)^k}{q_\theta(\Delta,X_i,X_{i+1})^k}.
\end{equation}
 An analytical development of $L_n^{(1)}(\theta)$ holds around $0$:
\begin{equation}\label{series}
L^{(1)}_n(\theta)=\sum_{k=0}^{+\infty}\theta^k L_n^{(k+1)}(0).
\end{equation}
Condition \eqref{eq:bound} implies that $|L_n^{(k)}(0)|\leq n$
and thus the series in \eqref{series} 
is absolutely convergent for $|\theta|<1$.

Introduce, for $k=1,2,\dots$, the sequence of functions
\begin{equation*}
h_k(x,y)=\left[\sgn(x+y)\exp\left(-(2/T)(x(x+y))^+\right)\right]^k.
\end{equation*}
We can then rewrite $L_n^{(k)}(0)$, for $k=1,2,\dots$, as: 
\begin{equation*}
L_n^{(k)}(0)=
(-1)^{k-1}\sum_{i=0}^{n-1}h_k(\sqrt{n}X_{i},\sqrt{n}(X_{i+1}-X_{i})).
\end{equation*}
We then see that the study of the limit behaviour of this type of sums,
presented in the next proposition, can be directly obtained from results 
obtained by J.~Jacod \cite{jacod}.
(For convenience, we present Jacod's results from \cite{jacod} in an Appendix).

\begin{proposition}
\label{prop-1}
Assume that $\theta=0$ in \eqref{eq:skb}, i.e.
let $X$ be a Brownian motion on $[0,T]$ departing from $x$, 
and let $\ell^x$ denote its local time at zero. 
\begin{itemize}
\item[\rm (a)]
Assume that $k=2,4,\dots$. 
Denote 
\begin{equation}\label{eq:ck}
\mu_k=-2\int_0^{\infty}\left[
1+{1\over 2k-1}\exp\left({2k(k-1)x^2\over(2k-1)^2}\right)\right]\Phi(-x)\,dx.
\end{equation}
Then 
\begin{equation}\label{eq:even}
\frac{L_n^{(k)}(0)}{n^{1/2}}
\xrightarrow[n\to\infty]{\text{prob.}}
\mu_k \ell^x_T.
\end{equation}
\item[\rm (b)]
Assume that $k=1,3,\dots$.
Denote 
\begin{equation*}
\mu_k=2\int_0^{\infty}\left[
1+{1\over 4k-1}\exp\left({4k(2k-1)x^2\over(4k-1)^2}\right)\right]\Phi(-x)\,dx.
\end{equation*}
Then, there exists a Brownian motion $W$ independent from $B$ such that  
\begin{equation}\label{eq:odd}
\frac{L_n^{(k)}(0)}{n^{1/4}}
\xrightarrow[n\to\infty]{\text{$\mathcal{F}$-stably}}{\mu}_{k}W(\ell^x_T). 
\end{equation}
\end{itemize}
\end{proposition}

\begin{remark}\label{remark:zero} Observe that on the event 
\begin{equation*}\label{eq:zero}
\{\omega\colon\inf_{0\leq t\leq T}X_t(\omega)>0\}
\end{equation*}
we have 
$\ell^x_T(\omega)=W\left(\ell^x_T(\omega)\right)=0$.
In this situation, as all the information about
the relevant parameter $\theta$ is produced when the process hits the level zero, no statistical inference can be carried out. Observe that in case $x=0$ we have 
\begin{equation*}
\PP\left( \left\{\omega\colon\inf_{0\leq t\leq T}X_t>0\right\}\right)=0.
\end{equation*}
\end{remark}

\begin{remark}
\label{rem:joint-conv}
Indeed, the results of J.~Jacod could be applied to multi-dimensional
statistics. This way, we obtain the joint $\mathcal{F}$-stable convergence 
of any vector $n^{-1/4}(L^{(1)}_n(0),\dotsc,L^{(2k+1)}_n(0))$
for any integer $k$,
and then the joint $\mathcal{F}$-stable convergence of
\begin{multline*}
(n^{-1/4}L^{(1)}_n(0),n^{-1/2}L^{(2)}_n(0),\dotsc, 
n^{-1/4}L^{(2k+1)}_n(0),n^{-1/2}L^{(2k+2)}_n(0))\\
\xrightarrow[n\to\infty]{\text{$\mathcal{F}$-stably}}
(\mu_1 W(\ell_T^x),\mu_2\ell^x_T,\dotsc,
\mu_{2k+1} W(\ell_T^x),\mu_{2+2}\ell^x_T).
\end{multline*}
\end{remark}

\begin{proof} 
We apply Theorem \ref{th1.1} in the Appendix. Observe that
\begin{equation*}
h_k(x,y)\leq \exp\left(-(x(x+y))^+\right)\leq \exp\left(|y|-|x\wedge x^2|\right).
\end{equation*}
We then have that \eqref{eq:r} holds with $a=1$, 
$\hat{h}(x)=\exp\left(-|x\wedge x^2|\right)$ and $r=0$, 
then it holds for any $r>0$. 
In consequence, by the aftermentioned Theorem,
the convergence in \eqref{eq:a} holds for $h=h_k$ with $k=2,4,\dots$. 
It rests to compute the constant in \eqref{eq:ch}.
We have
\begin{align}
c(h_k)&=\iint_{\RR^2}h_k(x,y)p(1,y)\,dx\,dy=
2\int_0^{\infty}\,dx\int_{-\infty}^{-x}p(1,y)\,dy\notag\\
   &\quad+{2\over\sqrt{2\pi}}
    \int_0^{\infty}\,dx\int_{-x}^{\infty}\exp\left(-\frac12y^2-2kxy-2kx^2\right)\,dy\notag\\
   &=2\int_0^{\infty}\Phi(-x)\,dx\notag\\
   &\quad+{2\over\sqrt{2\pi}}
    \int_0^{\infty}\,dx\exp\left(2k(k-1)x^2\right)\int_{-x}^{\infty}\exp\left(-\frac12(y+2kx)^2\right)\,dy\notag\\
   &=2\int_0^{\infty}\Phi(-x)\,dx+2
    \int_0^{\infty}\exp\left(2k(k-1)x^2\right)\Phi(-(2k-1)x)\,dx\notag\\
   &=2\int_0^{\infty}\Phi(-x)\,dx+{2\over 2k-1}\int_0^{\infty}\exp\left({2k(k-1)x^2\over(2k-1)^2}\right)\Phi(-x)\,dx.
\label{eq:chk}
\end{align}
Taking into account that $\mu_k=-c(h_k)$  
we conclude that \eqref{eq:even} holds with $\mu_k$ given in \eqref{eq:ck}.

To prove (b) we rely on Theorem \ref{th1.2} in the appendix. 
Observe then that $c(h_k)=0$ for odd $k$ due to the property
\begin{equation*}
h_k(-x,-y)=-h_k(x,y)\quad\text{for odd $k$}.
\end{equation*}
In view of the fact that \eqref{eq:r} holds for all $h_k$ with $r=4$, 
taking into account that  $(h_k)^2=h_{2k}$, 
we conclude that ${\mu}_k=\,c\left(h_{2k}\right)$. 
In view of the the computations in \eqref{eq:chk} with $2k$ instead of $k$ we conclude \eqref{eq:odd},
and the proof of the proposition. 
\end{proof}

\begin{remark}
On the event $A_n$, the discrete path has crossed the origin.
Hence, the continuous path did so at a random time $\tau$. Using the strong Markov 
property, this implies that one may consider a path starting
from $0$ for any time $t\geq \tau$. The local time of the Brownian 
motion is equal in distribution to the maximum of the Brownian motion.
Hence, on $A_n$ and $A$, $\ell_t^x>0$ for any time $t\geq\tau$.
\end{remark}

\begin{corollary}
\label{cor-1}
In the conditions of Proposition \ref{prop-1}, for $k=0,1,2,\dotsc$, we have:
\begin{equation}
{L^{(2k+1)}(0)\over L^{(2)}(0)}
\xrightarrow[n\to\infty]{\text{prob.}}
0,\qquad
{L^{(2k+2)}(0)\over L^{(2)}(0)}
\xrightarrow[n\to\infty]{\text{prob.}}
{\mu_{2k+2}\over\mu_2}\indicator_A,\label{eq:odd-even-limit}\\
\end{equation}
\begin{equation}
n^{1/4}{L^{(2k+1)}(0)\over L^{(2)}(0)}\mid A_n
\xrightarrow[n\to\infty]{\text{$\mathcal{F}$-stably}}
{\mu_{2k+1}W(\ell^x_T)\over\mu_2\ell^x_T}\mid A,\label{eq:dist}
\end{equation}
In particular, as $\mu_2=-\mu_1$, we obtain
\begin{equation}
-n^{1/4}{L^{(1)}(0)\over L^{(2)}(0)}\mid A_n
\xrightarrow[n\to\infty]{\text{$\mathcal{F}$-stably}}
{W(\ell^x_T)\over\ell^x_T}\mid A,\label{eq:alpha}
\end{equation}
\end{corollary}

\begin{proof} We begin by the second part in \eqref{eq:odd-even-limit}. As $\ell^x_T>0$ on the set $A$, we have 
\begin{equation*}
{L^{(2k+2)}(0)\over L^{(2)}(0)}\indicator_A
={n^{-1/2}L^{(2k+2)}(0)\over n^{-1/2}L^{(2)}(0)}\indicator_A
\xrightarrow[n\to\infty]{\text{prob.}}
{\mu_{2k+2}\ell^x_T\over\mu_2\ell^x_T}\indicator_A={\mu_{2k+2}\over\mu_2}\indicator_A.
\end{equation*}
Assume now that $\omega\in A^c$. We have $a(\omega)=\inf_{0\leq t\leq T}X_t(\omega)>0$ a.s. on $A^c$. Now
\begin{align*}
|L^{(2k+2)}_n(0)|&\leq (2k+1)!\sup_{1\leq i\leq n}e^{\left(-4kn/T\right)X_iX_{i+1}}|L_n^{(2)}(0)|\\
&\leq 
(2k+1)!e^{\left(-(4k-2)n/T\right)a(\omega)^2}|L_n^{(2)}(0)|,
\end{align*}
what gives the second part of \eqref{eq:odd-even-limit} on the set $A^c$. 
We postpone by now the proof of the first part of \eqref{eq:odd-even-limit}.

Let us then verify \eqref{eq:dist}.
We first prove that
\begin{equation*}
\left(
\frac{L_n^{(2)}(0)}{n^{1/2}},\frac{L^{(2k+1)}_n(0)}{n^{1/4}},\indicator_{A_n}
\right)
\xrightarrow[n\to\infty]{\text{$\mathcal{F}$-stably}}
\left(
\mu_2\ell^x_T,\mu_{2k+1}W(\ell^x_T),\indicator_A
\right).
\end{equation*}
This amounts to prove that, for $Z\geq 0$, $\mathcal{F}$-measurable and bounded,
and real $\lambda,\mu$ and $\nu$, we have
\begin{multline}\label{eq:preliminary}
\delta_n:=\EE
Z\exp
\left(i
\left\{
\lambda\frac{L_n^{(2)}(0)}{n^{1/2}}+\mu\frac{L^{(2k+1)}_n(0)}{n^{1/4}}+\nu\indicator_{A_n}
\right\}
\right)
\\
\xrightarrow[n\to\infty]{}
\EE
Z\exp
\left(
i\left\{
\lambda\mu_2\ell^x_T+\mu\mu_{2k+1}W(\ell^x_T)+\nu\indicator_A
\right\}
\right)=:\delta.
\end{multline}
We know that
\begin{equation*}
\left(\frac{L_n^{(2)}(0)}{n^{1/2}},\indicator_{A_n}\right)
\xrightarrow[n\to\infty]{\text{prob.}}
\left(\mu_2\ell^x_T,\indicator_{A}\right),
\end{equation*}
as $\indicator_{A_n}\to\indicator_{A}$ a.s. We then have
\begin{multline*}
|\delta_n-\delta|\leq
\left|
\EE
Z\exp
\left(i
\left\{
\lambda\frac{L_n^{(2)}(0)}{n^{1/2}}+\mu\frac{L^{(2k+1)}_n(0)}{n^{1/4}}+\nu\indicator_{A_n}
\right\}
\right)
\right.
\\
\left.
\qquad\qquad
-
\EE
Z\exp
\left(
i
\left\{
\lambda\mu_2\ell^x_T+\mu\frac{L^{(2k+1)}_n(0)}{n^{1/4}}+\nu\indicator_A
\right\}
\right)
\right|
\\
\hspace*{-1cm}
+
\left|
\EE
Z\exp
\left(
i\left\{
\lambda\mu_2\ell^x_T+\mu\frac{L^{(2k+1)}_n(0)}{n^{1/4}}  +\nu\indicator_A
\right\}
\right)
\right.\\
\qquad\qquad\qquad
\left.
-
\EE
Z\exp
\left(
i\left\{
\lambda\mu_2\ell^x_T+\mu\mu_{2k+1}W(\ell^x_T)+\nu\indicator_A
\right\}
\right)
\right|
\\
\leq
\EE
Z
\left|
\exp
\left(i
\left\{
\lambda\frac{L_n^{(2)}(0)}{n^{1/2}}+\nu\indicator_{A_n}
\right\}
\right)
-
\exp
\left(i
\left\{
\lambda\mu_2\ell^x_T+\nu\indicator_{A}
\right\}
\right)
\right|
\\
+
\EE
Z
\left|
\exp
\left(i
\mu\frac{L_n^{(2k+1)}(0)}{n^{1/4}}
\right)
-
\exp
\left(i
\mu\mu_{2k+1}W(\ell^x_T)
\right)
\right|
\to 0,
\end{multline*}
concluding the proof of \eqref{eq:preliminary}. The proof of
\eqref{eq:dist} follows with the help of the continuous and bounded
function 
$f_{K}(t)=
t\indicator_{\{|t|\leq K\}}+K\indicator_{\{t>K\}}-K\indicator_{\{t<-K\}}$. We have
\begin{equation*}
f_K\left(n^{1/4}{L^{(2k+1)}(0)\over L^{(2)}(0)}\indicator_{A_n}\right)
\xrightarrow[n\to\infty]{\text{$\mathcal{F}$-stably}}
f_K\left(
{\mu_{2k+1}W(\ell^x_T)\over\mu_2\ell^x_T}\indicator_A
\right)
\end{equation*}
for all $K>0$, and, as the limit is bounded in probability, we obtain \eqref{eq:dist}.

In what respects the first part of \eqref{eq:odd-even-limit} the computation on the set $A^c$ 
is similar to the previous one.
In the set $A$, we have
\begin{equation*}
{L^{(2k+1)}(0)\over L^{(2)}(0)}\indicator_A
=n^{-1/4}\left({n^{-1/4}L^{(2k+1)}(0)\over n^{-1/2}L^{(2)}(0)}\indicator_A\right)
\xrightarrow[n\to\infty]{\text{prob.}}
0,
\end{equation*}
as the expression within brackets has weak limit.
\end{proof}


\subsection{A simple estimator}

The MLE is the point $\theta_n$ at
which $\theta\mapsto L_n(\theta)$ 
reaches its maximum, i.e. $\theta_n$ is the (unique in our case) root of the equation
\begin{equation*}
L^{(1)}_n(\theta_n)=0.
\end{equation*}

Let us set 
\begin{equation}\label{alpha}
\alpha_n=-n^{1/4}\frac{L_n^{(1)}(0)}{L_n^{(2)}(0)}.
\end{equation}
From Corollary~\ref{cor-1} formula \eqref{eq:alpha}, 
$\alpha_n$ is known to converge $\mathcal{F}$-stably as $n\to\infty$.

Below, we will see that $n^{1/4}\theta_n$ and $n^{1/4}\alpha_n$ 
have the same limit, which yields Theorem~\ref{thm:1}.

The value $\alpha_n/n^{1/4}$, which is pretty simple to compute
from the data, specially in contrast to $\theta_n$ that requires a numerical solver to be computed,
can be used as an estimator of the skewness parameter.

Let us also remark that $\alpha_n$ is chosen so that the
first two terms in the Taylor series \eqref{series} of $L_n^{(1)}(\alpha_n/n^{1/4})$
cancel out. 


\subsection{Asymptotic development of the MLE}

We then prove a theorem and a theorem and a proposition which enclose Theorem~\ref{thm:1}.

\begin{theorem}
\label{thm:2} 
For any integer $p\geq 0$, there exists a vector
$(d_n^{(1)},\dotsc,d_n^{(p+1)})$ of random variables
given the recursive relation $d_n^{(1)}=1$ and 
\begin{equation}
\label{eq-10}
d_n^{(m+1)}=-
\sum_{k=1}^{m+1}\frac{L_n^{(k+1)}(0)}{L_n^{(2)}(0)}\sum_{\substack{1\leq i_1,\dotsc,i_k\leq m\\
i_1+\dotsb+i_k=m+1}}
d_n^{(i_1)}\dotsb d_n^{(i_k)}
\end{equation}
that converges $\mathcal{F}$-stably conditioning to $A$ to a vector 
$(d^{(1)},\dotsc,d^{(p+1)})$ depending only on $\ell_T^x$ and
$W(\ell_T^x)$. Besides, for any $\epsilon>0$,
there exists some integer $n_0$ large enough and some $K$ such that 
\begin{equation*}
\PP\left[n^{\frac{p}{4}+\frac{1}{2}}
|\theta_n-\Theta_n|\geq K\right]\leq \epsilon
\text{ for any }n\geq n_0, 
\end{equation*}
where 
\begin{equation*}
\Theta_n=\frac{\alpha_n}{n^{1/4}}
+d_n^{(2)}\frac{\alpha_n^2}{n^{1/2}}
+\dotsb
+d_n^{(p+1)}\frac{\alpha_n^{p+1}}{n^{(p+1)/4}}.
\end{equation*}
In addition, $d_n^{(2)}$ converges
to $0$ and $n^{1/4}d_n^{(2)}$ is bounded.
\end{theorem}

We prove this theorem after the next proposition, which 
will be stated in the following framework:
Using the result of Remark~\ref{rem:joint-conv}, 
we consider the asymptotic behavior of 
the vector 
$$
\left(n^{-1/4}L^{(1)}_n(0),n^{-1/2}L^{(2)}_n(0),\dotsc,n^{-1/2}L^{(2k)}_n(0)\right)
$$
for some $k\geq 1$.
We may then consider a probability space
$(\widehat{\Omega},\widehat{\mathcal{F}},\widehat{\PP})$ such that
this sequence is equal in distribution to a sequence 
converging almost surely to $(\mu_1 W(\ell_T^x),\mu_2 \ell_T^x,\dotsc,\mu_{2k}\ell_T^x)$.
We now consider some point in this probability space such that 
$\ell_T>0$. If the starting point is $0$, then the event
$\{\ell_T>0\}$ is of full measure. 

\begin{proposition}
\label{prop:2}
On the probability space 
$(\widehat{\Omega},\widehat{\mathcal{F}},\widehat{\PP})$ above, 
the random sequences $d^{(i)}_n$ given by \eqref{eq-10}
are convergent and bounded in $n$. Besides, 
for $m=1,2,3,\dotsc$,
\begin{equation*}
\theta_n=\frac{\alpha_n}{n^{1/4}}
+d_n^{(2)}\frac{\alpha_n^2}{n^{2/4}}
+d_n^{(3)}\frac{\alpha_n^3}{n^{3/4}}
+\dotsb
+d_n^{(m)}\frac{\alpha_n^m}{n^{m/4}}
+\grandO\left(\frac{1}{n^{(m+1)/4}}\right)
\end{equation*}
almost surely in the event $\{\ell_T^x>0\}$. 
\end{proposition}

Let us start by a simple lemma to get a control over
the finite Taylor expansion of $L^{(1)}_n(\theta)$.

\begin{lemma}
\label{lem-2}
For any $\theta$ and any integer $m\geq 1$, we have that 
for a random constant $C$ such that 
\begin{equation}
\label{eq-6}
\left|L^{(1)}_n(\theta)-\sum_{k=0}^{m}L_n^{(k+1)}(0)\theta^k\right|
\leq \sup_{|\xi|\leq |\theta|} |L_n^{(m+2)}(\xi)|\cdot|\theta|^{m+1}
\leq Cn^{1/2}\frac{|\theta|^{m+1}}{(1-|\theta|)^{m+2}}.
\end{equation}
\end{lemma}

\begin{proof}
With \eqref{eq:bound} and for $\theta\in(-1,1)$, 
\begin{equation}
\label{eq-8}
1-|\theta|
\leq 
\frac{q_\theta(\Delta,X_i,X_{i+1})}{q_0(\Delta(X_i,X_{i+1}))}\leq 1+|\theta|.
\end{equation}
with \eqref{derk}, since $L_n^{(k)}(0)/n^{1/2}$ converges in probability
(either to $\mu_k\ell_T^x$ or to $0$ depending if $k$ is even or odd),  
there exists a random constant $C$ 
\begin{equation*}
|L^{(k)}_n(\theta)|\leq \frac{n^{1/2}C}{(1-|\theta|)^k}.
\end{equation*}
Hence
\begin{equation*}
\left|L^{(1)}_n(\theta)-\sum_{k=0}^{m}L_n^{(k+1)}(0)\theta^k\right|
\leq \sup_{|\xi|\leq |\theta|} |L_n^{(m+2)}(\xi)|\cdot|\theta|^{m+1}.
\end{equation*}
With \eqref{eq:even} and \eqref{eq:odd}, this gives \eqref{eq-6}
because $L_n^{(k)}(0)/n^{1/2}$ is bounded in $n$.
\end{proof}

\begin{lemma}
\label{lem-3}
For $n$ large enough, the function $L_n^{(1)}(\theta)$ is invertible.
Besides, the function $(L_n^{(1)}(\theta))^{-1}$ is Lipschitz in $\theta$
with a constant $8/n^{1/2}\mu_2\ell^x_T$ on the event $A_n$.
\end{lemma}

\begin{proof}
With \eqref{eq-8},
\begin{equation*}
-L_n^{(2)}(\theta)\geq \left(\sum_{i=1}^{n-1}
\frac{p(\Delta,|X_i|+|X_{i+1}|)^2}{q_0(\Delta,X_i,X_{i+1})^2}
\right)\frac{1}{(1+|\theta|)^2}
\geq \frac{L_n^{(2)}(0)}{4}.
\end{equation*}
Since $L^{(2)}_n(0)<0$ for $n$ large enough 
as $n^{-1/2}L^{(2)}_n(\theta)$ converges in probability to
some negative random variable (see \eqref{eq:even}), 
we get that $L_n^{(1)}(\theta)$ is one-to-one.
With the formula $\partial_\theta (L^{(1)}(\theta))^{-1}=1/L_n^{(2)}(L_n^{(1)}(\theta))$,
$(L_n^{(1)}(\theta))^{-1}$ is Lipschitz in $\theta$ with 
constant $4/L_n^{(2)}(0)$. 
\end{proof}

The idea of the proof is then the following: 
We construct an of estimator~$\Theta_n$ such that 
for some constant $C$ and $p\geq 0$, 
\begin{equation*}
\sup_{n\in\NN} n^{p/4} |L_n^{(1)}(\Theta_n)|\leq C.
\end{equation*}
Since $L_n^{(1)}(\theta_n)=0$, 
\begin{multline}
\label{eq-9}
|\Theta_n-\theta_n|=|(L_n^{(1)})^{-1}(L_n^{(1)}(\Theta_n))
-(L_n^{(1)})^{-1}(L_n^{(1)}(\theta_n))|\\
\leq \frac{8}{n^{1/2}\mu_2\ell^x_T}|L_n^{(1)}(\Theta_n)|
\leq \frac{8C}{n^{p/4+1/2}\mu_2\ell^x_T}.
\end{multline}

\begin{proof}[Proof of Proposition~\ref{prop:2}]
For the sake of simplicity, let us set $q=n^{1/4}$. 

Set $\Theta_n=\alpha_n q+\beta_n q^2+\gamma_n q^3+\xi_n q^4$
for some $\beta_n$, $\gamma_n$ and~$\xi_n$ to be carefully chosen.
Here, we consider only the first terms in the development 
of $\theta_n$. It is easily to convince oneself that this 
method may be applied to any order and that the involved
terms $\beta_n,\gamma_n,\xi_n,\dots$ 
may be computed recursively
and gives rise to \eqref{eq-10}.

With \eqref{eq-6} and $m=4$, there exists a constant $C$ such that
\begin{multline}
\label{eq-7}
|L^{(1)}(\Theta_n)
-L_n^{(1)}(0)-L_n^{(2)}(0)\Theta_n\\
-L_n^{(3)}(0)\Theta_n^2
-L_n^{(4)}(0)\Theta_n^3
-L_n^{(5)}(0)\Theta_n^4
|
\leq Cn^{1/2}\frac{|\Theta_n|^5}{
(1-|\Theta_n|)^6}.
\end{multline}
Remark that $L_n^{(1)}(0)-L_n^{(2)}(0)\alpha_n q=0$.
In order to get rid of the terms in $q^2$, set
\begin{equation*}
\beta_n=\frac{-L_n^{(3)}(0)}{L_n^{(2)}(0)}\alpha_n^2
\end{equation*}
Since $\alpha_n$ converges and $n^{1/4}L_n^{(3)}(0)/L_n^{(2)}(0)$
also converges stably, then $n^{1/4}\beta_n$ converges stably.
Also, $\beta_n$ converges to $0$.

In order to get rid of the terms in $q^3$, set
\begin{equation*}
\gamma_n=-\frac{L_n^{(4)}(0)}{L_n^{(2)}(0)}\alpha_n^3
-\frac{L_n^{(3)}(0)}{L_n^{(2)}(0)}\alpha_n\beta_n.
\end{equation*}
From Corollary~\ref{cor-1}, $\gamma_n$ converges stably since 
$\alpha_n$ and $\beta_n$ converges stably.

In order to get rid of the terms in $q^4$, set
\begin{equation*}
\xi_n=-2\frac{L^{(3)}_n(0)}{L^{(2)}_n(0)}(\alpha_n\gamma_n+\beta_n^2)
-4\frac{L^{(4)}_n(0)}{L^{(2)}_n(0)}\alpha_n^2\beta_n.
\end{equation*}
Again, $\xi_n$ converges thanks to Corollary~\ref{cor-1}.

Hence
\begin{equation*}
L_n^{(1)}(\Theta_n)=\sum_{r=5}^{20} q^r B_n^{(r)}+R_n(\Theta_n),
\end{equation*}
where $R_n(\Theta_n)\leq n^{1/2}|\Theta_n|^{5}/(1-|\Theta_n|)^6$
and $B_n^{(r)}$ are terms that depend 
linearly on $L^{(k)}_n(0)$ and on the power 
of the $\alpha_n$, $\beta_n$, $\gamma_n$ and $\xi_n$.
Since the $L^{(k)}_n(0)/n^{1/2}$ are bounded, 
we obtain that the $n^{3/4}B_n^{(r)}$ are bounded.

In addition, $n^{1/4}\Theta_n$ is bounded in $n$, so 
that $n^{3/4}R_n(\Theta_n)$ is bounded in $n$.
With \eqref{eq-9}, 
this proves that for some constant $K$, 
\begin{equation*}
|\Theta_n-\theta_n|\leq \frac{K}{n^{5/4}}.
\end{equation*}
This result may be generalized to any order.
Finally, let us note that 
$\beta_n=\alpha_n^2 d_n^{(2)}$ with $d_n^{(2)}
=-L_n^{(3)}(0)/L_n^{(2)}(0)$. With \eqref{eq:odd-even-limit}, 
$d_n^{(2)}$ converges in probability to $0$
and $n^{1/4}d_n^{(2)}$ is bounded in $n$.
$\gamma_n=\alpha_n^3 d_n^{(3)}$ and
$\xi_n=\alpha_n^4 d_n^{(4)}$ where 
$d_n^{(3)}$ and $d_n^{(4)}$ are bounded in $n$.
\end{proof}

\begin{proof}[Proof of Theorem~\ref{thm:1}]
Let us consider the event $\{\ell_T^x>0\}$. It corresponds
to the event $A$ as the local time of the Brownian motion is positive
just after having hit $0$.
Since on this event, the local time has a density (see Lemma~\ref{lem-1} below)
which is derived from the one of the first hitting time 
of a point $x$, for each $\epsilon>0$, one may find
a set $\Omega(\epsilon)$ as well as some values $0<a'<b'$
and $c'$ such that 
$\omega\in\Omega(\epsilon)$ implies that 
$\ell_T^x\in(a',b')$ and $|W(\ell_T(x))/\ell_T^x|\leq c'$
and 
\begin{equation*}
\PP[\Omega(\epsilon)|\{\ell_T^x>0\}]\geq 1-\epsilon/2.
\end{equation*}

From the joint convergence of the $n^{1/4}L_n^{(2k+1)}(0)$ 
to $\mu_{2k+1} W(\ell_T^x)$ and the joint convergence
of the $n^{1/2}L_n^{(2k)}(0)$ to $\mu_{2k}\ell_T^x$, we get 
for any $\epsilon>0$, there exists $0<k<a'$ and $K>b'$ 
as well as a measurable set
$\Omega'(\epsilon,n)\subset \Omega(\epsilon)$ such that 
\begin{align*}
L_n^{(k)}(0)&\leq K\sqrt{n}\text{ on }\Omega'(\epsilon,n)\\
L_n^{(2)}(0)&\geq k\sqrt{n}\text{ on }\Omega'(\epsilon,n)\\
\text{ and }
\forall n\geq n_0,&\ \PP[\Omega'(n,\epsilon)|\{\ell_T^x>0\}]\geq 1-\epsilon.
\end{align*}

In the proof of Proposition~\ref{prop:2}, we constructed some
estimator $\Theta_n$ such that for some $p\geq 0$, $n^p L_n^{(1)}(\Theta_n)$
is bounded by some constant depending the upper bounds
of the $n^{1/2}L_n^{(k)}(0)$. Besides, we use the Lipschitz
constant of $(L_n^{(1)})^{-1}$ which depends on the lower
bound of $n^{1/2}L_n^{(2)}(0)$. Thus, on $\Omega'(\epsilon,n)$, 
we obtain that $|\theta_n-\Theta_n|\leq C/n^{p+1/2}$, where 
$C$ depends only on $K$ and $k$, assuming that $n\geq n_0$.
This means that 
\begin{equation*}
n^{p+1/2}|\theta_n-\Theta_n|\leq C.
\end{equation*}
Thus, for any $\epsilon>0$, there exists $n_0$ large
enough such that 
\begin{equation*}
\forall n\geq n_0,\ \PP[n^{p+1/2}|\theta_n-\Theta_n|\geq C]\leq \epsilon.
\end{equation*}
which yields the result.
\end{proof}

\subsection{The contrast function}

In order to study the maximum likelihood, it is also possible 
to consider the contrast function
\begin{equation*}
Z_n(\theta)=\frac{\exp(L_n(\theta))}{\exp(L_n(0))}.
\end{equation*}
Using the asymptotic development of $L_n(u)$ around $0$,
we get that 
\begin{equation*}
\log Z_n(\theta)=\theta L_n^{(1)}(0)
+\frac{\theta^2}{2}L_n^{(2)}(0)
+\grandO(\theta^3).
\end{equation*}
Thus, with the result of Proposition~\ref{prop-1} and taking into account that $\mu_2=-{\mu}_1$,
we see that 
\begin{equation}\label{eq:lamn}
\log Z_n(\theta/n^{1/4})\xrightarrow[n\to\infty]{\text{$\mathcal{F}$-stably}}
\mu_1\left(\theta W(\ell_T^x)-{\theta^2\over 2}\ell_T^x\right).
\end{equation}
From this convergence we can intuitively check our result in \eqref{eq:theta}, based in the theory of
convergence of statistical experiments and the LAMN property in \eqref{eq:lamn}. The theory states
(under certain stringent conditions that we do not verify)
that the maximum likelihood estimator of the pre-limit experiments converges stably to the maximum
likelihood estimator of the limit experiment \cite{ih}. It is direct, differentiating with respect to $\theta$
in the r.h.s. of \eqref{eq:lamn}, to obtain, when $\ell_T^x>0$, that the MLE in the limit experiment
is $W(\ell_T^x)/\ell_T^x$. We then obtain \eqref{eq:theta} in the form
\begin{equation*}
n^{1/4}\theta_n
\xrightarrow[n\to\infty]{\text{$\mathcal{F}$-stably}}
W(\ell_T^x)/\ell_T^x.
\end{equation*}

\section{The limit distribution}\label{section:ld}

As $n^{1/4}\theta_n$ and $n^{1/4}\alpha_n$ 
converge to $\Upsilon=W(\ell^x_T)/\ell^x_T$, we give 
the main characteristics of this random variables.
To simplify the computations, we assume
that $x=0$ and $T=1$, so that we write $\ell_1=\ell^x_T$.

Indeed, this random variable is easy to simulate.

\begin{lemma}
\label{lem-1}
The distribution of $\Upsilon$ is symmetric.
Besides, its density is 
\begin{equation}
\label{eq-3}
f_\Upsilon(x)=\frac{\dd F_\Upsilon(x)}{\dd x}
=\int_0^{+\infty}\vd y\int_0^1 
\frac{\sqrt{y}}{2\pi\sqrt{t^3}}
\exp\left(\frac{-xy}{2}-\frac{y^2}{2t}\right)\vd t.
\end{equation}
and it is equal in distribution to 
\begin{equation}
\label{eq-4}
\Upsilon=\frac{G(H)}{H}
\text{ with }
H=\frac{1}{2}(U+\sqrt{V+U^2}),
\end{equation}
where $G(H)$, $U$ and $V$ are independent random variables, 
$G(H)$, $G(H)\sim\cN(0,H)$, $U\sim\cN(0,1)$ and 
$V\sim\exp(1/2)$.
\end{lemma}

\begin{proof}
It is well known that the local time $\ell_1$ at time $1$
is equal in distribution to the supremum of the Brownian motion 
$\sup_{r\in[0,1]} B_r$
on $[0,1]$. It follows that
\begin{equation*}
F_{\ell_1}(y)=\PP_0[\ell_1<y]=\PP_0[\sup_{r\in[0,1]} B_r<y]
=\PP_0[\tau_y>1],
\end{equation*}
where $\tau_y=\inf\{t>0\,|\,B_t=y\}$. The density $\upsilon(t;y)$
of $\tau_y$ is equal to 
\begin{equation*}
\upsilon(t;y)=\frac{1}{\sqrt{2\pi t^3}}\exp\left(-\frac{y^2}{2t}\right),
\end{equation*}
so that 
\begin{equation*}
F_{\ell_1}(y)=1-\int_0^{1}
\frac{1}{\sqrt{2\pi t^3}}\exp\left(-\frac{y^2}{2t}\right)\vd t,
\end{equation*}
and the density $f_{\ell_1}(y)$ of $\ell_1$ is then equal to 
\begin{equation*}
f_{\ell_1}(y)=\int_0^{1}
\frac{2y}{\sqrt{2\pi t^3}}\exp\left(-\frac{y^2}{2t}\right)\vd t.
\end{equation*}
Thus, conditioning with respect to the value of $\ell_1$,
\begin{equation*}
F_\Upsilon(x)=\PP[\Upsilon<x]
=\int_0^{+\infty} \PP[W(y)<xy]f(y)\vd y
\end{equation*}
and this leads to \eqref{eq-3}.

Expression \eqref{eq-4} follows from the equality
in distribution of $\ell_1$ and $\frac{1}{2}(U+\sqrt{V+U^2})$.
This expression has been used in order to simulate
the reflected Brownian motion \cites{lepingle93a,lepingle95a}.
\end{proof}

\begin{figure}[H]
\begin{center}
\includegraphics[angle=-90,scale=0.6]{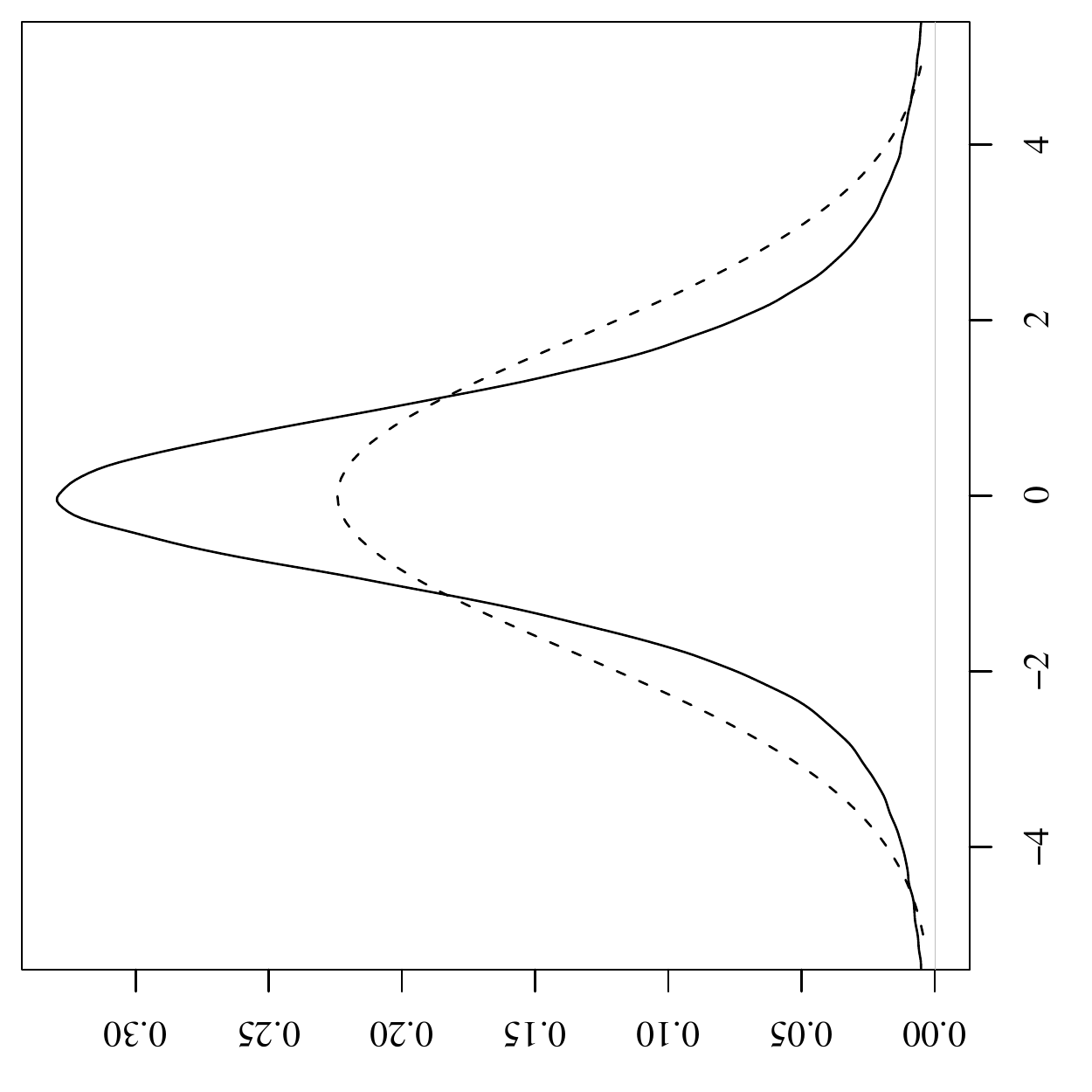}
\caption{\label{fig-1} Density of $\Upsilon$ (solid)
and density of the normal distribution with variance
$\Var(\Upsilon)$ (dashed).}
\end{center}
\end{figure}

The variance of $\Upsilon$ is $3.16$.
We see in Figure~\ref{fig-1} that the density of $\Upsilon$
is close to that of the normal distribution, yet narrower.

\section{Numerical tests and observations on the likelihood}

Numerical tests are easy to perform, as all the formulae are
easy to implement.

\subsection{On the coefficient $\alpha_n$}

Several tests can be performed on $\alpha_n=-n^{1/4}L^{(1)}_n(0)/L^{(2)}_n(0)$, mainly to see whether it is
reasonable to use it instead of the MLE $\theta_n$.

First, one can check that $\theta_n$ and $\alpha_n/n^{1/4}$ are pretty close, 
by setting $\theta_n=\argmax_{\theta\in(-1,1)} L_n(\theta)$
and computing it using a numerical procedure.
In Table~\ref{table-1}, one can check that the
error of $|\theta_n-\alpha_n/n^{1/4}|$ is of order $1/n^{3/4}$,
so that $\alpha_n/n^{1/4}$ is a pretty good approximation of $\theta_n$, 
and is much more faster to compute.

\begin{table}
\begin{center}
\begin{tabular}{r r c c c c}
\hline
$n$ & mean & $n^{1/2}\times$ mean & $n^{3/4}\times$ mean  & std dev & quant. $90\,\%$   \\ 
\hline
\np{100} &  \np{0.026} & \np{0.26} & \np{0.8} & \np{0.057} & \np{0.082} \\
\np{200} &  \np{0.028} & \np{0.40} & \np{1.5} & \np{0.083} & \np{0.057} \\
\np{500} &  \np{0.013} & \np{0.29} & \np{1.3}  &\np{0.055} & \np{0.026} \\
\np{1000} & \np{0.013} & \np{0.41} & \np{2.3} &\np{0.040} & \np{0.033} \\
\np{2000} & \np{0.006} & \np{0.26} & \np{1.8}  &\np{0.025} & \np{0.015} \\
\np{5000} & \np{0.006} & \np{0.42} & \np{3.5} &\np{0.041} & \np{0.006} \\
\np{10000} & \np{0.002} & \np{0.20} & \np{2.0} &\np{0.005} & \np{0.003} \\
\hline
\end{tabular}
\caption{\label{table-1} Statistics of $|\theta_n-\alpha_n/n^{1/4}|$
over $100$ paths.}
\end{center}
\end{table}

Second, one can check the variance of $\alpha_n$, as 
well as the adequacy of $\alpha_n$ with the distribution 
of $\Upsilon=W(\ell_1)/\ell_1$. For this, we have
used a set of \np{10000} simulations of $\Upsilon$, 
and we have renormalized $\Upsilon$ to get
the same variance as $\alpha_n$. 
Using a Kolmogorov-Smirnov test, we can see in Table~\ref{table-2} 
that even for a low value of $n$ (\textit{e.g.} $n=1000$), 
we get a good adequation with the distribution of $\Upsilon$.
Yet, for $n=100$, the distribution of $\alpha_n/n^{1/4}$
or $\theta_n$ (by keeping only the values in $(-1,1)$, which
means $88\%$ of the values of $\alpha_n$ with $n=100$ and $96\%$ for
$n=\np{1000}$)
is in fact close to the Gaussian distribution. 

\begin{figure}[hp!]
\begin{center}
\subfloat[$\theta_n$ and $\alpha_n$ against normal, $n=100$]{\includegraphics[angle=-90,scale=0.5]{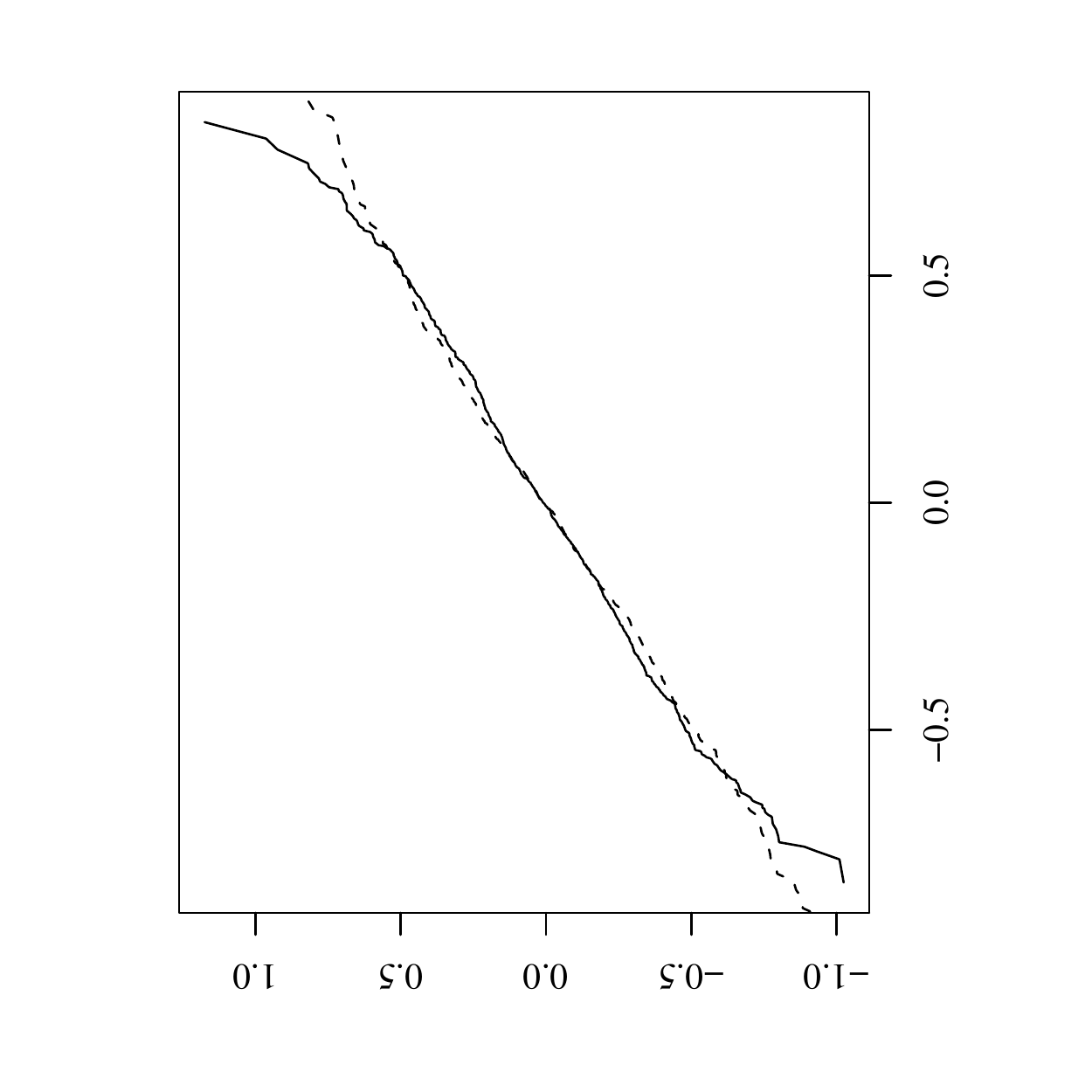}}
\subfloat[$\theta_n$ and $\alpha_n$ against $\Upsilon$, $n=100$]{\includegraphics[angle=-90,scale=0.5]{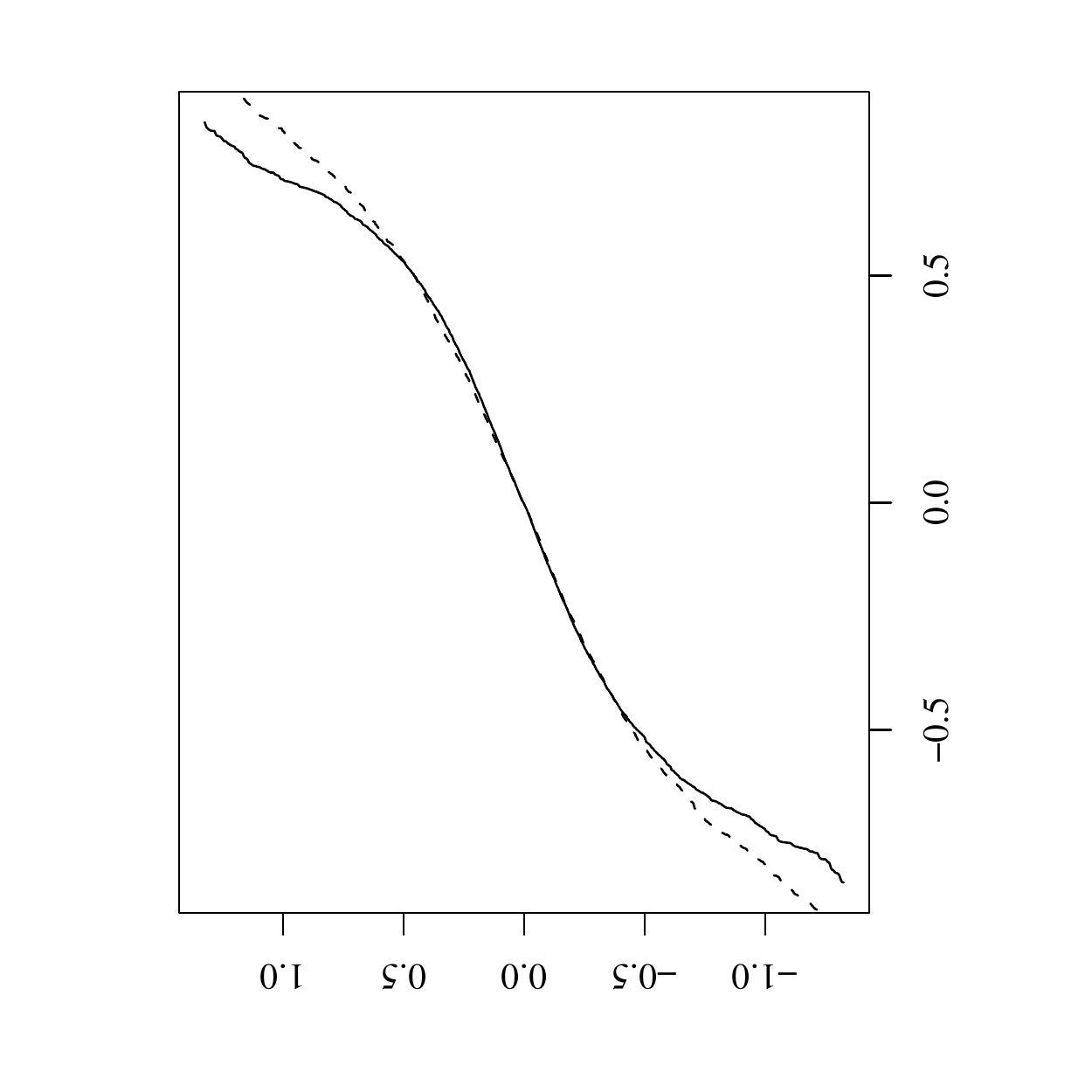}}

\subfloat[$\theta_n$ and $\alpha_n$ against normal, $n=1000$]{\includegraphics[angle=-90,scale=0.5]{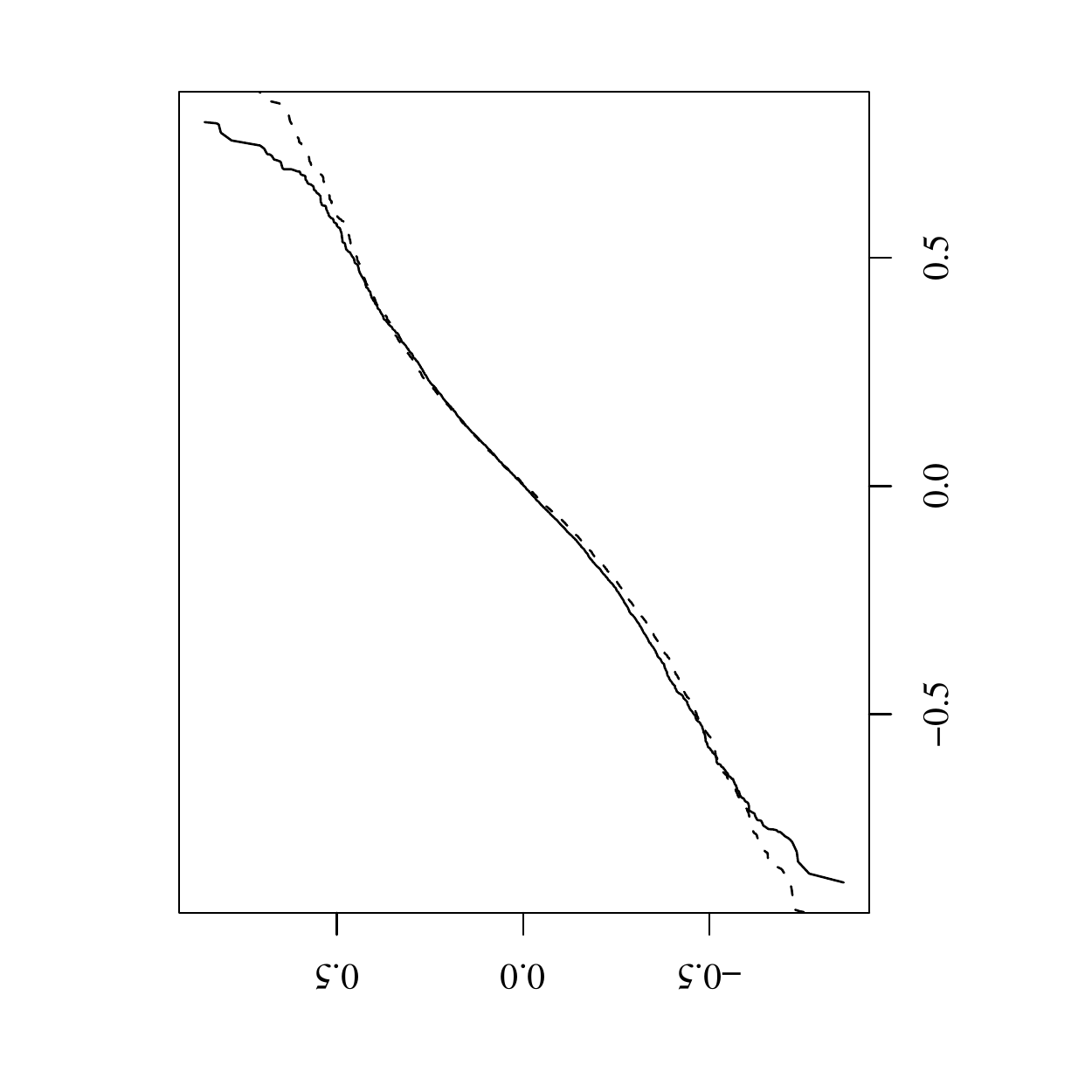}}
\subfloat[$\theta_n$ and $\alpha_n$ against $\Upsilon$, $n=1000$]{\includegraphics[angle=-90,scale=0.5]{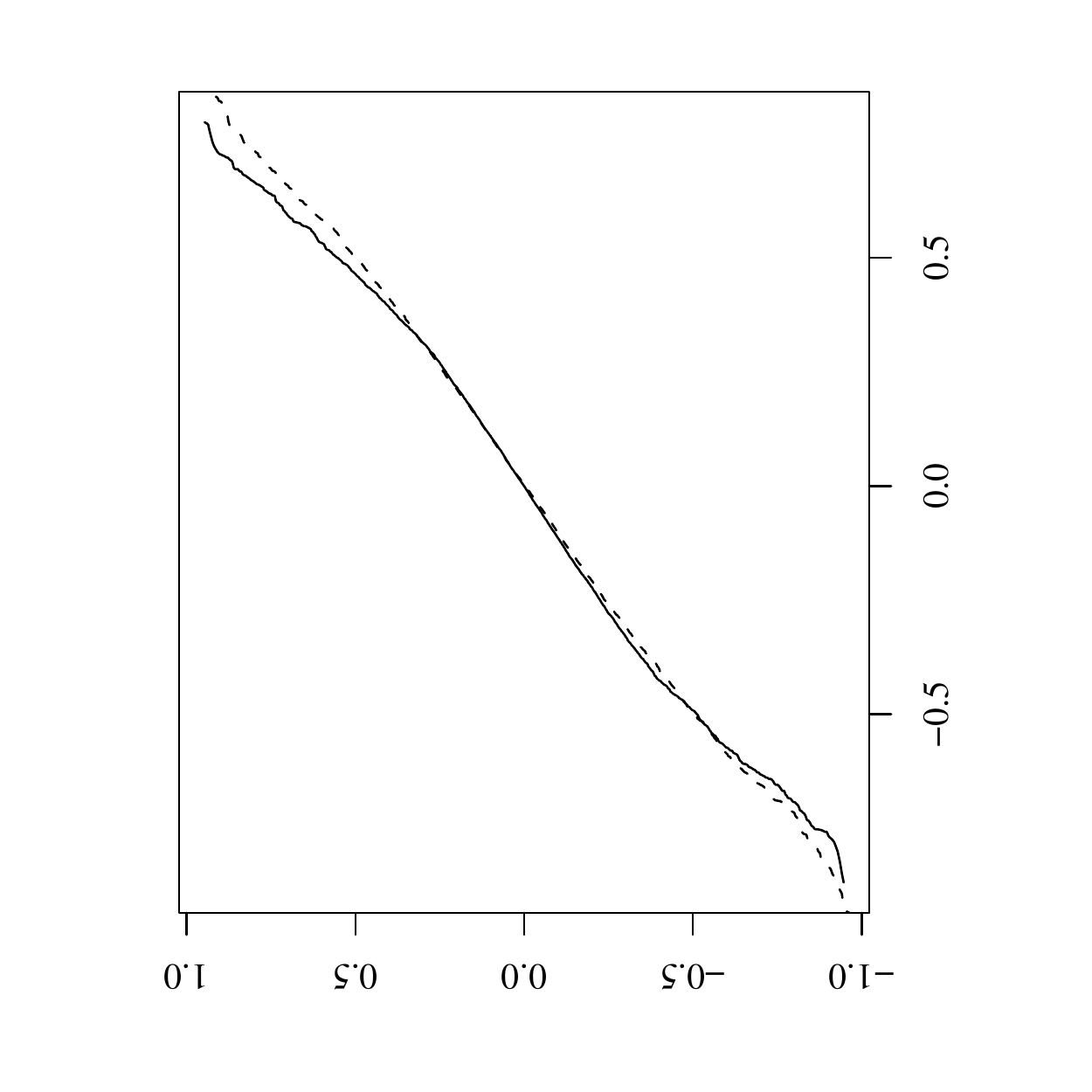}}
\caption{\label{fig-2} Quantile-Quantile
plot of $\theta_n$ (solid) and $\alpha_n/n^{1/4}$ (dashed)
against the normal distribution (left) and 
the distribution of $\Upsilon$ (right) with variance $\Var(\theta_n)$
for \np{10000} samples.}
\end{center}
\end{figure}

\begin{table}
\begin{center}
\begin{tabular}{|r| r r | r r|}
\hline
$n$ &   $D_n(\alpha_n,\Upsilon)$ & $p$-value  & $D_n(\alpha_n,G)$ & $p$-value\\
\hline
\np{100} & \np{0.020} & \np{0.30} & \np{0.082} & $<$\np{2e-16}\\
\np{250} & \np{0.017} & \np{0.005} & \np{0.079} & $<$\np{2e-16}\\
\np{500} & \np{0.024} & \np{1e-5} & \np{0.083} & $<$\np{2e-16}\\
\np{750} & \np{0.029} & \np{9e-8} & \np{0.087} & $<$\np{2e-16}\\
\np{1000} & \np{0.012} & \np{0.098} & \np{0.072} & $<$\np{2e-16}\\
\np{2500} & \np{0.023} & \np{4e-5} & \np{0.085} & $<$\np{2e-16}\\
\np{5000} & \np{0.021} & \np{2e-4} & \np{0.079} & $<$\np{2e-16} \\
\hline
\end{tabular}
\caption{\label{table-2} Kolmogorov-Smirnov test on $\alpha_n$
against $\Upsilon\sqrt{\Var(\alpha_n)/\Var(\Upsilon)}$
and the normal distribution $G$ with variance $\Var(\alpha_n)$
over \np{10000} paths.}
\end{center}
\end{table}

However, the variance of $\alpha_n$ is dependent on $n$
and is not stable with $n$.

In addition, for small values of $n$, there are
some values of $\alpha_n$ such that $\alpha_n/n^{1/4}$
is outside $[-1,1]$.


\subsection{On the order of convergence}

\begin{figure}[H]
\begin{center}
\includegraphics[angle=-90,scale=0.4]{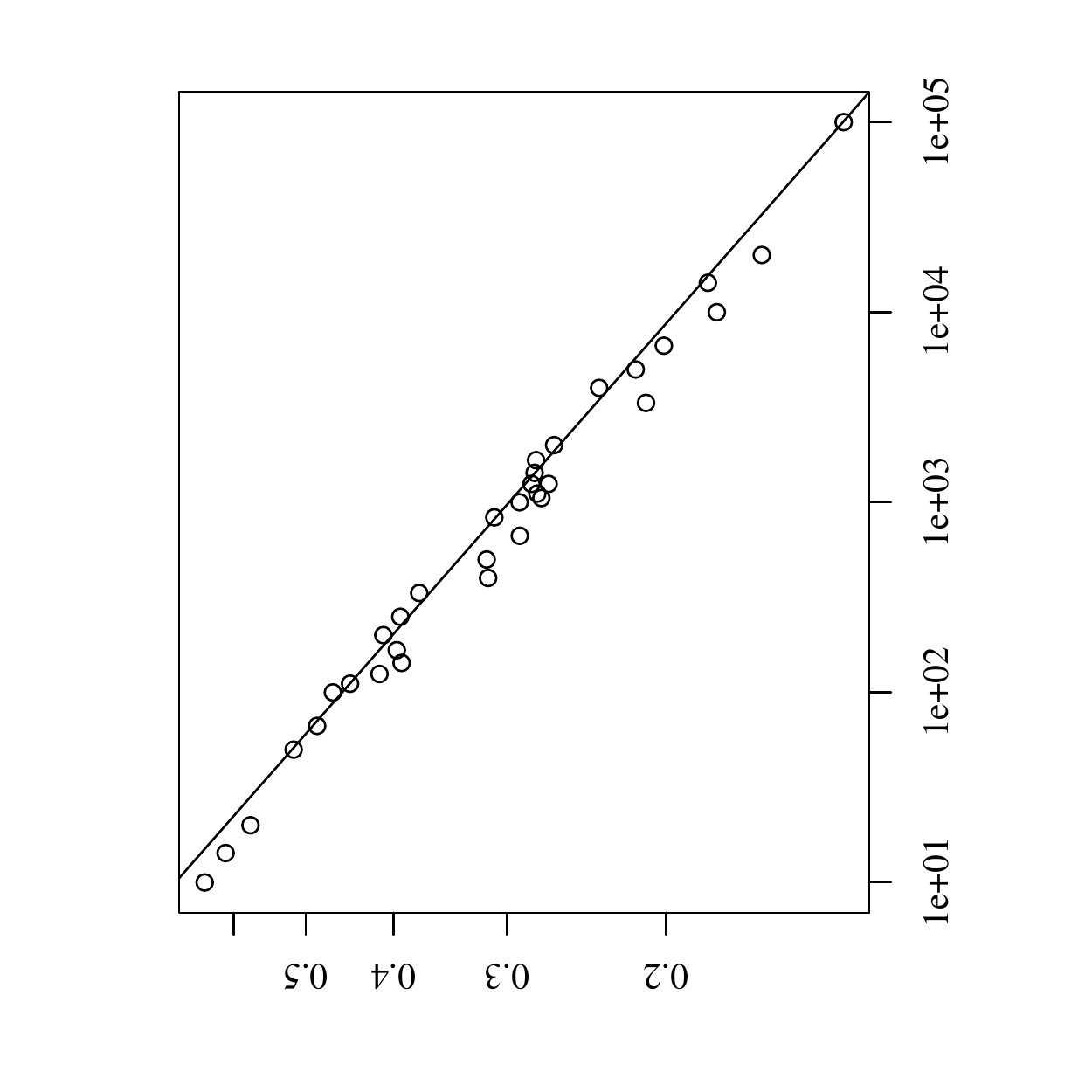}
\caption{\label{fig-3} Logarithmic regression on the standard 
deviation of $\theta_n$}
\end{center}
\end{figure}

One could wonder if the rate of convergence of $\theta_n$ is really
of order  $-1/4$.
Numerical simulations show that the rate of convergence, 
for $n$ in the range 100 to \np{100000} is of order $\delta$ 
with $\delta\approx -0.18$, which is smaller than $-0.25$.
This value is found using a regression on the logarithm
standard deviation of 500 samples of~$\theta_n$ (See Figure~\ref{fig-3}).

Indeed, one can note that the variance of $\alpha_n$
also depends on $n$, and in the range from 
\np{50} to \np{600000}, a numerical study of \np{10000}
samples of $\alpha_n$ shows that $\Var(\alpha_n)$
seems to be equal to $Cn^{\beta}$ with $\beta\approx 0.08$.
This has to be taken into account in order to design
some test of hypotheses.

\begin{figure}[H]
\begin{center}
\hbox to 0pt{\hfil
\includegraphics[scale=0.5,angle=-90]{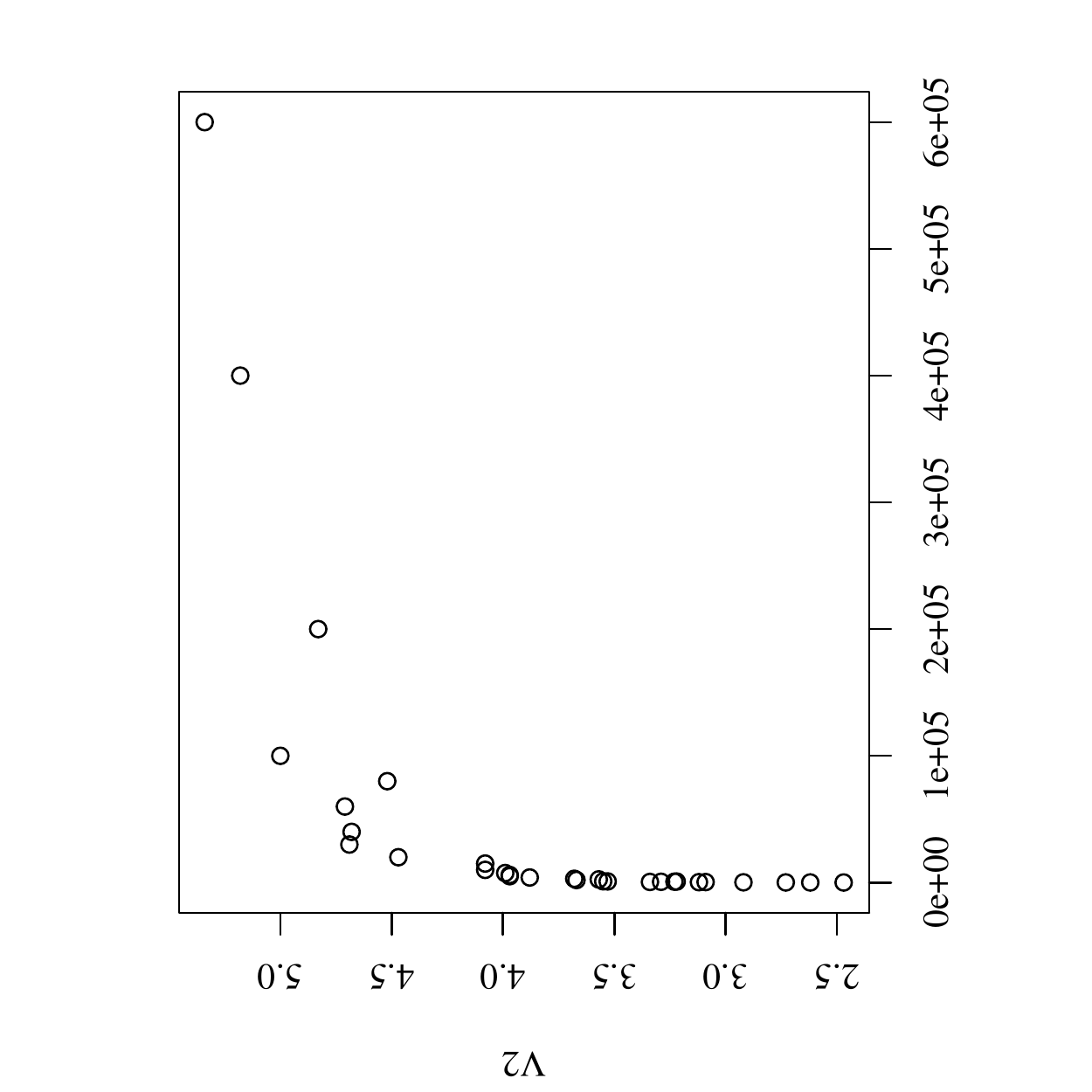}
\hfil
\includegraphics[scale=0.5,angle=-90]{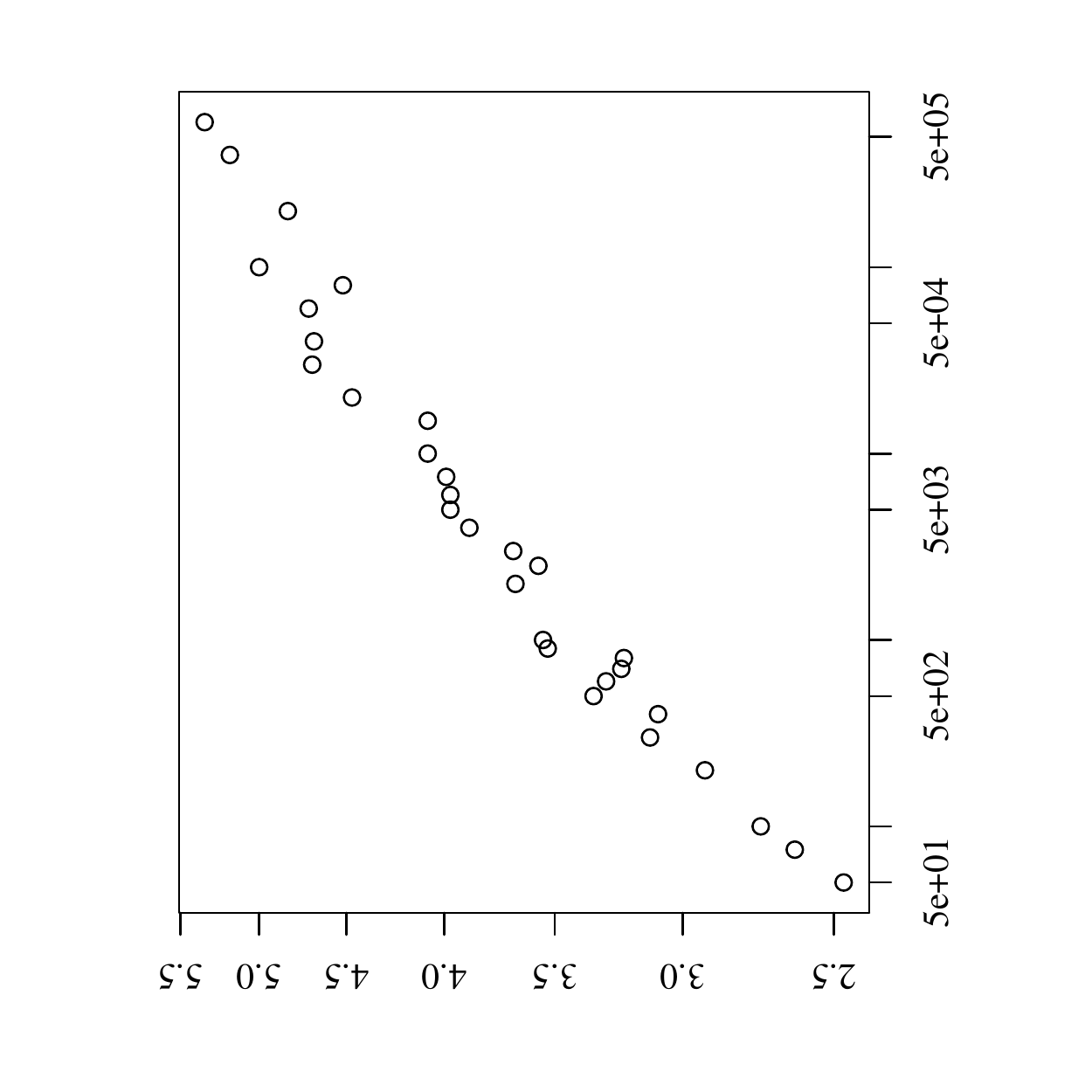}\hfil}

\caption{\label{fig-4} Variance of $\alpha_n$ as a function 
of $n$ in the linear and logarithmic scale.}
\end{center}
\end{figure}

\subsection{A hypothesis test}

It is then possible to develop a hypothesis test
of $\theta=0$ against $\theta\not=0$. For this, 
let us compute
\begin{equation*}
\PP\left[|\theta_n|\geq \frac{K}{n^{1/4}}\right]
\approx 
\PP\left[\left|\Upsilon\right|\geq \frac{K}{cn^{1/4}}\right].
\end{equation*}

\begin{figure}[H]
\begin{center}
\includegraphics[angle=-90,scale=0.6]{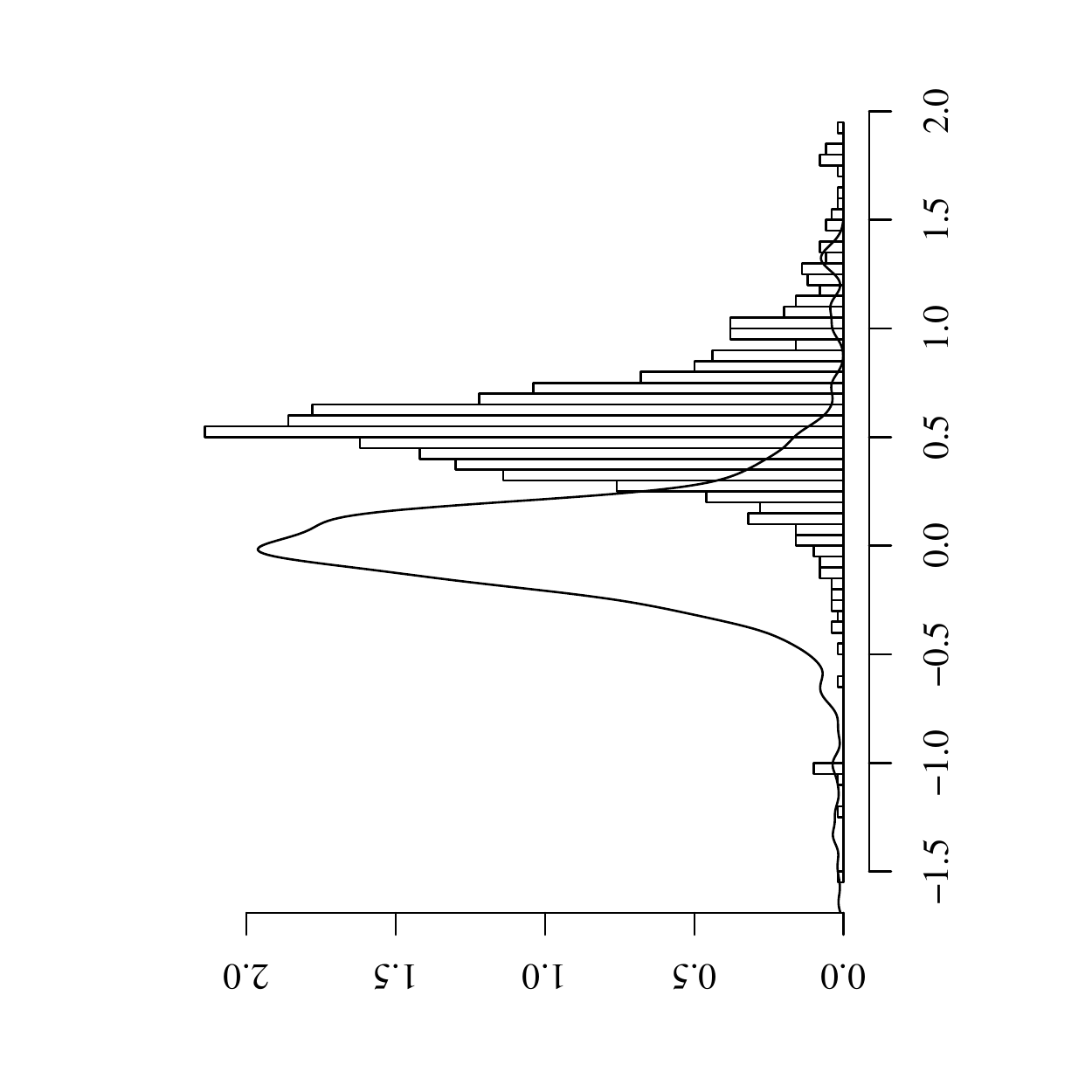}

\caption{\label{fig-5} Histogram of the density of $\alpha_n/n^{1/4}$
for $n=\np{1000}$ realizations of the SBm with $\theta=0.5$
against the an approximation of the density of $\alpha_n/n^{1/4}$
for the Brownian motion.}
\end{center}
\end{figure}

Of course, the second type error cannot be computed, 
and we do not know the asymptotic behavior of $L_n(\theta)$
when $\theta\not=0$. However, it is rather easy 
to perform simulation and thus to get some 
numerical information about the MLE $\theta_n$
and $\alpha_n$. For example, we see in Figure~\ref{fig-5}
an approximation of the density of $\alpha_n/n^{1/4}$
for $\theta=0.5$ compared to an approximation of the 
density of $\alpha_n/n^{1/4}$ for the Brownian motion
with $n=\np{1000}$. We can note that 
the histogram of $\alpha_n/n^{1/4}$ has its peak 
on $0.5$.

\section{An example of application: diffusion of species}

As endowed in the introduction, 
the SBm is a fundamental tool when one has to model a permeable barrier. 
In addition, 
it appears when one writes down the processes generated by diffusion 
equations with discontinuous coefficients in a one dimensional media: 
this issue is presented in the survey article \cite{lejay} with references to the articles where the SBm arised
and covering various fields, such as ecology, finance,
astrophysics, geophysics, ...

We present here a possible application to ecology
of our hypothesis test, which can be surely applied to other fields.

\subsection{Has a boundary between two habitats an effect?}

Diffusions are commonly used in ecology to explain the spread
of a specie, at the level of individual cells 
(See for example the book \cite{berg})
or the level of an animal in a wild environment.

Several authors have proposed the use of biased diffusions 
to model the behavior of a specie at the boundary  between two habitats
\cites{cantrell99a,ovaskainen03a}), when the species diffuse
with different species at speed in each habitat. 

Now, consider a situation where the dispersion of a
specie in two different habitats is well modelled by 
a diffusion process, and that the measurement of 
the diffusion coefficient give the same value.
Does it means that the boundary has no effect on the
displacement of the individuals?

Let us apply this in a one-dimensional world, 
where one habitat is $[0,+\infty)$ and the other
is $(-\infty,0]$. We assume that we may track
the position of an individual, whose displacement
in each of the habitat is given by $x+\sigma B_t$.

Then, we may apply our hypothesis test to determine
whether or not the position 
shall be modelled by 
\begin{equation*}
(\mathrm{H}_0)\quad X_t=x+\sigma B_t
\end{equation*}
or by 
\begin{equation*}
(\mathrm{H}_1)\quad X_t=x+\sigma B_t+\theta \ell_t^0(X).
\end{equation*}
Under Hypothesis $(\mathrm{H}_0)$, the boundary
has no effect and is not seen. 
Under Hypothesis $(\mathrm{H}_1)$, the individual
is more likely to go in one of the two habitat, depending
on the sign of $\theta$.

\subsection{What is the underlying operator?}

Now, let us consider that we have a measurement 
of the diffusion coefficients that gives 
two different values $a_+$ on $\RR_+$ and $a_-$ on $\RR_-$.

One may then wonder which differential operator shall 
be used to model the diffusive behavior. 
For $a=a_+\mathbf{1}_{[0,+\infty)} +a_-\mathbf{1}_{(-\infty,0)}$,
is it 
\begin{equation*}
L=\frac{1}{2}\nabla(a\nabla\cdot)\text{ or }A=\frac{1}{2}a\triangle?
\end{equation*}
On $(0,+\infty)$ and $(-\infty,0)$, there is no difference
between these two operators, which means that the local 
dynamic of the particle/individual is not affected by
the choice of $L$ and $A$. However, the difference arises
at $0$: the process $X$ generated by $L$ is solution to 
\begin{equation*}
X_t=x+\int_0^t \sqrt{a(X_s)}\vd B_s+\frac{a_+-a_-}{a_++a_-}\ell_t^0(X)
\end{equation*}
while the process $Y$ generated by $A$ is solution to 
\begin{equation*}
Y_t=x+\int_0^t \sqrt{a(X_s)}\vd B_s,
\end{equation*}
for a Brownian motion $B$ (See for example \cites{lejay-martinez04a,lejay}). 
From the analytical point of view: the domain $\Dom(A)$ of $A$ contains
the functions of class $\cC^2(\RR)$ which are bounded with 
bounded, first and second order derivatives. The domain $\Dom(L)$
of $L$ contains functions of class $\cC^2(\RR\setminus\{0\})$
with bounded first and second order derivatives which
are furthermore continuous at $0$, and such that 
$a_+\nabla f(0+)=a_-\nabla f(0-)$. This condition is called 
the \emph{flux condition}. In many physical situations, 
it is assumed that the flux $a\nabla u$ is continuous
and this is why divergence-form operators of type $L$ arise.

\begin{remark}
Both $L$ and $A$ can be embedded in a single class of operators
of type $\frac{\rho}{2}\nabla (a\nabla\cdot)$. If $\rho$ 
and $a$ are constant on $(0,+\infty)$ and $(-\infty,0)$, 
then we may use the following characterization: let us consider 
\begin{equation*}
C=\frac{1}{2}\nabla(a\nabla\cdot)\text{ with }a=a_+\mathbf{1}_{[0,+\infty)}
+a_-\mathbf{1}_{(-\infty,0)}
\end{equation*}
and 
\begin{equation*}
\Dom(C)=\left\{
f\in\cC^2(\RR\setminus\{0\})
\,\vrule\,
\begin{aligned}
&f,f',f''\text{ are bounded on $\RR\setminus\{0\}$}\\
&f(0-)=f(0+)\\
&(1+\lambda)f'(0+)=(1-\lambda)f'(0-),\ \lambda\in(-1,1)
\end{aligned}
\right\}.
\end{equation*}
This class of operators is then specified by three parameters, 
$a_+>0$, $a_->0$ and~$\lambda\in(-1,1)$. The operator $A$ corresponds
to $\lambda=0$, while $L$ corresponds to 
$\lambda=(a_+-a_-)/(a_++a_-)$.
\end{remark}

For $\Phi(x)=\int_0^x \vd x/\sqrt{a(x)}$, 
$\widehat{X}=\Phi(X)$ is solution to the SDE \cites{lejay-martinez04a,lejay}
\begin{equation*}
\widehat{X}_t=\Phi(x)+B_t+\frac{\sqrt{a_+}-\sqrt{a_-}}{\sqrt{a_+}+\sqrt{a_-}}\ell_t^0(\widehat{X}),
\end{equation*}
while $\widehat{Y}=\Phi(Y)$ is solution to the SDE
\begin{equation*}
\widehat{Y}_t=\Phi(x)+B_t+\frac{\sqrt{a_-}-\sqrt{a_+}}{\sqrt{a_+}+\sqrt{a_-}}\ell_t^0(\widehat{Y}).
\end{equation*}
We then see that both $\widehat{X}$ and $\widehat{Y}$ are
Skew Brownian motions, but the coefficients in front of their local
time have opposite signs.

Even if we have not studied the asymptotic behavior of the
MLE for the SBm with skewness parameter different from $0$, 
numerical experiments back the following hypotheses test:
\begin{enumerate}
\item Given an observed $X$, estimate the diffusion coefficient for 
the process on each side of $0$.
\item Apply the function $\Phi$ to the observed process.
\item Compute the MLE $\theta_n$ of the Skewness parameter.
If $a_+>a_-$ (resp. $a_+<a_-$) and $\theta_n>0$ then 
decided that the infinitesimal generator of~$X$ is $L$ (resp. $A$).
Otherwise, decide that it is $A$ (resp. $L$).
\end{enumerate}

\section{Conclusion}

In this article, we have studied the behavior
of the maximum likelihood for the Skew Brownian
motion when the parameter to estimate is $0$.

In particular, we have shown that the rate
of convergence of the estimator~$\theta_n$
is $n^{1/4}$ and not $n^{1/2}$ as in the
classical case. 
This should not be surprising: 
indeed, away from~$0$, the Skew Brownian motion
behaves like a Brownian motion, and only its
dynamic close to~$0$ allows one to see
the difference between a Skew Brownian 
motion with a parameter $\theta\not=0$ and
a Brownian motion. It is also not surprising
that the local time enters in the limit distribution.

The case $\theta\not=0$ remains open. One needs to prove results similar to the
one of J.~Jacod \cite{jacod}, when the Brownian motion is replaced by the Skew
Brownian motion (its distribution with respect to the Wiener measure is
singular).  Of course, one cannot expect the limit law to be symmetric. Yet, it
is pretty easy to simulate the Skew Brownian motion and to estimate the maximum
likelihood, so that numerical studies are easy to perform.

\section{Appendix}
In this Appendix  we provide the theorems given in \cite{jacod} used for the proofs of the main results in Section 2.
We slightly change the notation and present the results in the particular cases that are relevant to us
in the present work.

Denote by $X=\{X_t\colon 0\leq t\leq T\}$ 
a Brownian motion on a probability space 
$(\Omega,\mathcal{F},\PP)$. Introduce a Borel function $h\colon\RR^2\to\RR$ such that
there exist $a\in\RR$ and $\hat{h}\colon\RR\to\RR$ such that 
\begin{equation}\label{eq:r}
h(x,y)\leq e^{a|y|}\hat{h}(x)\quad \text{and}\quad \int|x|^r|\hat{h}(x)|\,dx<\infty.
\end{equation}
\begin{theorem}[{\bf Theorem 1.1 p. 508 in \cite{jacod}}]\label{th1.1}
Consider $h$ as above, satisfying \eqref{eq:r} with $r=0$. Then
\begin{equation}\label{eq:a}
\frac1{n^{1/2}}\sum_{i=0}^{n-1}h(\sqrt{n}X_{i/n},\sqrt{n}(X_{(i+1)/n}-X_{i/n}))
\xrightarrow[n\to\infty]
{\text{prob.}}c(h)\ell^x_T,
\end{equation}
where
\begin{equation}\label{eq:ch}
c(h)=\iint_{\RR^2}h(x,y)p(1,y)\,dx\,dy,
\end{equation}
and $\ell^x$ denotes the local time of $X$ at level zero.
\end{theorem}
\begin{remark} It must be noticed that the convergence in \eqref{eq:a}, as stated in \cite{jacod},
is stronger, in the sense that both terms in \eqref{eq:a} are processes (i.e. depend on $t$) and
the convergence is locally uniformly in time, in probability.
Recall that a sequence $(Z^n)_{n\ge 1}$ of processes is said to \emph{converges 
locally uniformly in time, in probability}, to a limiting processes $Z$ if for 
any $t \in \mathbb{R}^+$ the sequence $\sup_{s \le t}|Z^n_s - Z_s|$ goes to $0$ 
in probability.
\end{remark}
\begin{theorem}[{\bf Theorem 1.2 p. 511 in \cite{jacod}}]\label{th1.2}
Consider $h$ as above, satisfying \eqref{eq:r} with some $r>3$, and assume that $c(h)=0$ (see \eqref{eq:ch}). 
Then
\begin{equation}\label{eq:b}
\frac1{n^{1/4}}\sum_{i=0}^{n-1}h(\sqrt{n}X_{i/n},\sqrt{n}(X_{(i+1)/n}-X_{i/n}))
\xrightarrow[n\to\infty]
{\text{$\mathcal{F}$-stable in dist.}}c\left(h^2\right) W(\ell^x_T),
\end{equation}
where $W=\{W_t\colon t\geq 0\}$ is a Brownian motion independent of $X$, and
$\ell^x$ is the local time of $X$ at level zero. 
The constant $c\left(h^2\right)$ is given in \eqref{eq:ch} for the
function $h^2$.
\end{theorem}

\begin{remark} As in the previous remark, the Theorem stated in \cite{jacod} is stronger,
now in the sense that both terms in \eqref{eq:b} are processes, and the processes converge
stably in distribution in the Skorokhod space of c\`adl\`ag functions.
\end{remark}



\begin{bibdiv}

\begin{biblist}

\bib{bardou}{article}{
    author={Bardou, O.},
    author={Martinez, M.},  
    title={Personal communication},
}
\bib{berg}{book}{
    author={Berg, H.C.},
    title={Random walks in Biology},
    year={1993},
    publisher={Princeton University Press},
}
\bib{cantrell99a}{article}{
    author={Cantrell, R.S.},
    author={Cosner, C.},
    title={Diffusion models for population dynamics incorporating individual
  behavior at boundaries: Applications to refuge design},
    journal={Theoretical Population Biology}, 
    volume={55},
    number={2},
    pages={189--207},
    year={1999},
}
\bib{florens1}{article}{
    author={Florens-Zmirou, D.},
    title={On estimating the diffusion coefficient from discrete observations},
    journal={J. Appl. Probab.},
    volume={30},
    number={4},
    pages={790--804},
    year={1993},
}

\bib{florens2}{article}{
    author={Florens-Zmirou, D.},
    title={Statistics on crossings of discretized diffusions and local time},
    journal={Stochastic Process. Appl.},
    year={1993},
    volume={39},
    pages={139--151}, 
    number={1},
}

\bib{harrison}{article}{
    author={Harrison, J.M.},
    author={Shepp, L.A.},
    title={On Skew Brownian motion},
    journal={Ann. Probab.},
    volume={9},
    number={2},
    year={1981},    
    pages={309--313},
}
\bib{ih}{book}{
author={Ibragimov I.A.},
author={Has'minskii R.Z.}, 
title={Statistical Estimation Asymptotic Theory},
publisher={Springer},
address={New York},
year={1981},
}

\bib{jacod}{article}{
    author={Jacod, J.},
    title={Rates of convergence to the local time of a diffusion},
    journal={Ann. Inst. H. Poincaré Probab. Statist.},
    year={1998},
    volume={34},
    number={4},
    pages={505-544},
}

\bib{jacod2}{article}{
    author={Jacod, J.},
    title={Parametric inference for discretely observed non-ergodic diffusions},
    journal={Bernoulli},
    volume={12},
    number={3},
    year={2006},
    pages={383--401},
}
\bib{jacod-shiryaev:1987}{book}{
author={Jacod, J.},
author={Shiryaev, A. N.},
title={Limit theorems for stochastic processes.},
publisher={Springer-Verlag},
year={1987},
}

\bib{K}{book}{
    author={Kutoyants, Yu.A.},
    title={Parameter Estimation for Stochastic Processes},
    publisher={Heldermann},
    address={Berlin},
    year={1984},
}
\bib{lecam-yang}{book}{
   author={Le Cam, L.},
   author={Yang, G. L.},
   title={Asymptotics in statistics},
   series={Springer Series in Statistics},
   edition={2},
   publisher={Springer-Verlag},
   place={New York},
   date={2000},
}
\bib{lejay}{article}{
    author={Lejay, A.},
    title={On the constructions of the Skew Brownian motion},
    journal={Probab. Surv.},
    year={2006},
    volume={3},
    pages={413--466},
}

\bib{lejay-martinez04a}{article}{
    author={Lejay, A.},
    author={Martinez, M.},
    title={A scheme for simulating one-dimensional diffusion processes with
  discontinuous coefficients},
    journal={Ann. Appl. Probab.},
    volume={16},
    number={1},
    pages={107--139}, 
    year={2006},
}

\bib{lepingle93a}{article}{
   author = {L{\'e}pingle, D.},
    title = {Un sch\'ema d'{E}uler pour \'equations diff\'erentielles
             stochastiques r\'efl\'echies},
  journal = {C. R. Acad. Sci. Paris S\'er. I Math.},
   volume = {316},
     year = {1993},
   number = {6},
    pages = {601--605},
}

\bib{lepingle95a}{article}{
    author = {L{\'e}pingle, D.},
     title = {Euler scheme for reflected stochastic differential equations},
      note = {Probabilit\'es num\'eriques (Paris, 1992)},
   journal = {Math. Comput. Simulation},
    volume = {38},
      year = {1995},
    number = {1-3},
     pages = {119--126},
}

\bib{lipster}{book}{
    author={Lipster, R.S.},
    author={Shiryaev, A.N.},
    title={Statistics of random processes. II. Applications},
    publisher={Springer}, 
    address={Berlin},
    year={2001},
}
\bib{martinez04a}{thesis}{
	author = {Martinez, M.},
	title = {Interpr{\'e}tations probabilistes d'op{\'e}rateurs
sous forme divergence et analyse de m{\'e}thodes num{\'e}riques
associ{\'e}es},
	year = {2004},
	type = {Ph.D. thesis},
	school = {Universit{\'e} de Provence / INRIA Sophia-Antipolis},
}

\bib{ovaskainen03a}{article}{
    author={Ovaskainen, O.},
    author={Cornell, S.~J.}, 
    title={Biased movement at a boundary and conditional occupancy times for
  diffusion processes},
    journal={J. Appl. Probab.},
    volume={40},
    number={3},
    pages={557--580}, 
    year={2003},
}    
    

\end{biblist}
\end{bibdiv}
\end{document}